\documentclass{amsart}

\usepackage[all]{xy}
\usepackage{graphics}
\usepackage{hyperref}
\usepackage{amssymb}
\usepackage{amsmath}
\usepackage{amsthm}

\newtheorem{theo}{Theorem}[section]
\newtheorem{prop}{Proposition}[section]
\newtheorem{lemma}{Lemma}[section]
\newtheorem{cor}{Corollary}[section]

\newcommand{\varthm}[1]{\par\medskip\noindent\textbf{#1.}\ }
\newcommand{\tend}{\longrightarrow}

\newcommand{\adh}{\overline}
\newcommand{\ins}{\overset{\circ}}

\newcommand{\hbl}{\mathfrak{h}}
\newcommand{\bbl}{\mathfrak{b}}
\newcommand{\phibl}{\mathfrak{p}}

\newlength{\mywd}\newlength{\myht}\newlength{\myaux}
\newcommand{\crc}{{\ensuremath{\scriptstyle\circ}}}
\newcommand{\insb}[1]{%
 \settowidth{\mywd}{$#1$}%
 \settoheight{\myht}{$#1$}%
 \settowidth{\myaux}{\crc}
 \setlength{\myaux}{-0.5\myaux}
 \addtolength{\myaux}{0.5\mywd}%
 \raisebox{1.05\myht}[\myht][0pt]{%
  \makebox[0pt][l]{%
  \hspace{\myaux}%
  \crc%
  }%
 }%
 {{#1}}%
}

\newcommand{\bbC}{\mathbb{C}}
\newcommand{\bbD}{\mathbb{D}}

\newcommand{\bbH}{\mathbb{H}}

\newcommand{\bbN}{\mathbb{N}}

\newcommand{\bbR}{\mathbb{R}}
\newcommand{\bbS}{\mathbb{S}}
\newcommand{\bbT}{\mathbb{T}}

\newcommand{\bbZ}{\mathbb{Z}}

\newcommand{\mcB}{\mathcal{B}}

\newcommand{\mcP}{\mathcal{P}}

\begin{document}

\title{Ghys-like models providing trick for a class of simple maps}

\author{Arnaud Ch\'eritat}
\date{}

\maketitle
\noindent\textsc{\small Laboratoire \'Emile Picard, Universit\'e Paul Sabatier,
118~route de Narbonne, 31062~Toulouse Cedex, France}

\begin{abstract}
For quadratic polynomials with an indifferent fixed point with bounded type rotation number (they have a Siegel disk), much of what is known of their Julia set comes from the study of a quasiconformal model. The model is build from a Blaschke fraction, that we call a pre-model, and that is given by a formula. We give here a geometric construction of pre-model maps, that extends to some cases where no formula is known. More precisely, we are able to make this work for a class of entire maps, very specific but nonetheless spanning uncountably many equivalence classes (thus with probably no hope for a formula), and also in the case of the Lavaurs maps that arise in the parabolic implosion of quadratic polynomials.
\end{abstract}


\bf
\tableofcontents
\mdseries


\section{Introduction}

\subsection{Setting}

Let $\theta \in \bbR$ be irrational and $\rho = e^{i 2\pi\theta}$. Let
$f$ be a holomorphic function defined in a neighborhood of $0$, fixing
it with multiplier $\rho$, i.e.\ $f(0)=0$ and $f'(0)=\rho$. The
tangent map of $f$ at $0$ is an aperiodic rotation. Is $f$
linearizable, i.e.\ conjugated to it's linear part in a neighborhood
of $0$~? Not necessarily: G.\ A.\ Pfeiffer provided counter examples
in 1919, and Cr\'emer proved in 1927 that for ``bad values'' of
$\theta$, rational maps are not linearizable. On the other hand,
Siegel proved in 1942 that there are values of $\theta$ such that
every map $f$ as above is linearizable. More precisely, he proved this
for Diophantine numbers, i.e.\ number $\theta$ satisfying the
following condition: $\exists \sigma \geq 2$ and $C>0$, $\forall p\in
\bbZ$ and $q \in \bbN^*$, $|\theta - p/q | \geq C/q^\sigma$. Since,
the set of irrational $\theta$ such that every map $f$ is linearizable
has been completely characterized (by Brjuno \cite{Br} and Yoccoz
\cite{Yo}). It is a remarkable fact that it coincides with the set of
irrational $\theta$ such that the polynomial $P(z)= \rho z+z^2$ is
linearizable \cite{Yo}.

If $f$ is linearizable, there is a maximal domain where it is
conjugated to the rotation. This domain is a topological disk, called
the \textbf{Siegel disk}. Let us focus our attention on the family
$P(z) = \rho z + z^2$. Even in this case, the nature of the boundary
of the Siegel disk is not completely known. Obviously, it is
connected. If locally connected, then it is a Jordan curve (Rodin,
\cite{Rod}). On the other hand, the possibility that it be an
indecomposable continuum (there be no covering by two closed connected
proper subsets) for some quadratic polynomial is not completely ruled
out (see Rogers, \cite{Rog}). It might even happen that the boundary
of the Siegel disk is all the Julia set~! 

The first examples of Siegel disks with Jordan curve boundary were
given by Herman \cite{He1}. They are provided by a
\emph{quasiconformal surgery} procedure due to Ghys \cite{Gh}. The
surgery is done on a Blaschke fraction, i.e.\ a rational map which
preserves the unit circle, from which is constructed a \emph{model}
map. This model map, defined on the Riemann sphere, is quasiregular
and preserves a Beltrami form. It is conjugated to a rotation on the
closed disk, with rotation number $\theta$, but not beyond. One then
straightens the Beltrami form. This conjugates the model to some
holomorphic function on the Riemann sphere, and there is not much
choice: one proves this must be the polynomial $P$ (up to M\"obius
conjugacy). As an immediate consequence, the straightening maps the
unit disk to the Siegel disk $\Delta$ of $P$ and the unit circle to
the boundary of $\Delta$. As a corollary, the boundary of $\Delta$ is
a Jordan curve.

In this model, there is no critical points on the unit circle, thus
the boundaries of the Siegel disks of these examples do not contain
the critical point of $P$. Douady remarked in a Bourbaki
seminar~\cite{DBou} in 1987 that one should be able to apply the same
procedure with a Blaschke fraction $B$ having a critical point on the
unit circle $\bbS^1$, provided $B$ was quasisymmetrically conjugated
to the rotation on $\bbS^1$. Herman \cite{He2}, using inequalities by
\'Swi\c{a}tek \cite{Sw}, then proved in 1987 that this is the case if
and only if $\theta$ is of bounded type, i.e.\ Diophantine of
exponent~$\sigma=2$, or equivalently that the coefficients of the
continued fraction expansion of $\theta$ are bounded. This provided
examples where the critical point belongs to the boundary of the
Siegel disk.

This model proved useful for understanding the combinatorics of the
polynomial $P$ on it's Julia set $J$, and other properties. Working on
the model, Petersen proved that for bounded type rotation numbers, $J$
is locally connected~\cite{Pe1}. Then Petersen and Lyubich proved that
$J$ has Lebesgue measure zero.  (Petersen and Zakeri extended these
results to a wider class of irrational numbers, using Guy~David maps
instead of quasiconformal maps. In particular, this new class has full
Lebesgue measure in $\bbR$, contrarily to the class of bounded type
irrationals. See~\cite{PeZa}.) For bounded type rotation numbers,
McMullen~\cite{McMu} proved that the Julia set is shallow, which
implies it's Hausdorff dimension is $< 2$.

It is worth mentioning that these surgeries prove the linearizability
of $P$ at the fixed point $0$, without using Siegel's linearization
theorem.

Working with other models enables to transfer these results to other
holomorphic maps. See for instance~\cite{YaZa}, \cite{Ge}.

Recently, Graczyk and \'Swi\c{a}tek proved the following strong
result: for any map $f$ fixing a point with bounded type multiplier,
the Siegel disk must either have a critical point in it's boundary, or not be
compactly contained in the domain of definition of $f$, or both.
In particular, for any rational map, the first case holds.

\subsection{Purpose}

Usually, model maps are introduced by a formula. 
We provide here a geometric construction.
Does any function with a bounded type indifferent fixed point have a
model~? We will see that the geometric construction gives a positive
answer for functions having only one critical value and some form of
rigidity.

These models are tools for studying the corresponding Julia sets. Here
we define the tools and only give obvious corollaries. Further statements
(local connexity, measure, Hausdorff dimension) require more work.

We also warn the reader that for simplicity, we will not give all the details
in the case of Lavaurs maps. They can be found in~\cite{C}.


\section{The classical surgery}

We will call $\bbS^1$ the unit circle in the complex plane $\bbC$, and $\bbS^2$ the Riemann sphere.
Let $B$ be the rational map given by the following formula
\[
  B(z) = z^2\frac{z-3}{1-3z}
\]
This is a Blaschke fraction: $B(\bbS^1) = \bbS^1$. The restriction of
$B$ to $\bbS^1$ is an orientation preserving homeomorphism. It is not a
diffeomorphism, since the derivative cancels at $z=1$. The local
degree of this critical point is $3$. Let us call $B$ the \textbf{pre-model}.

Let $R_\tau$ be the rotation of center $0$ and rotation number
$\tau$: $R_\tau(z) = e^{i 2\pi \tau} z$.
We introduce the family $B_\tau = R_\tau \circ B$. According to the
theory of rotation number (see for instance~\cite{dMvS}),
for any irrational $\theta \in \bbR/\bbZ$, there is one and only one
value of $\tau \in \bbR/\bbZ$ such that $B_\tau$ has rotation number
$\theta$ on $\bbS^1$. According to the Herman-\'Swi\c{a}tek
theorem\footnote{see section~\ref{sec_premodeltomodel} for
statements}, the restriction of $B_\tau$ to $\bbS^1$ is
quasisymmetrically conjugated to $R_\theta$ if and only if $\theta$
has bounded type. If this is the case, then according to the
Ahlfors-Beurling theorem\footnote{idem},
the conjugacy extends to a quasiconformal self-homeomorphism $\phi$ of
the unit disk. The \textbf{model} $\widetilde{B}$ is defined by
\[
 \begin{array}{rclcl}
   \widetilde{B}(z) & = & B(z) &  & \text{if } z\in
   \bbS^2 \setminus \bbD \\
             & & \phi^{-1} \circ R_\theta \circ
    \phi (z) & & \text{if } z\in \bbD
 \end{array}
\]
This map is quasiregular. It is a ramified covering of degree $2$,
with critical points $\infty$ and $1$.
Let $\mu_0$ be the pull-back by $\phi$ of
the null Beltrami form on $\bbD$, and $\mu$ be the following
$\widetilde{B}$-invariant Beltrami form: for $z\in\bbS^2$, either the orbit
of $z$ under $B$ never gets in $\bbD$, in which case we set
$\mu(z)=0$, or there is a first visit $B^n(z) \in \bbD$ with
$n\in\bbN$, in which case we set $\mu(z)$ to be the value at $z$ of
the pull-back of $\mu_0$ by $B^n$.

Now we straighten $\mu_0$: there exists a quasiconformal homeomorphism $S : \bbS^2 \to \bbS^2$ that pushes $\mu_0$ to $0$ ($S_* \mu_0 = 0$), unique if we require that $S$ sends $\phi^{-1}(0)$ to $0$, $\infty$ to $\infty$, and $1$ to $-\rho/2$ where $\rho = e^{i 2\pi\theta}$. The function
\[
 F = S \circ \widetilde{B} \circ S^{-1} : \bbS^2 \to \bbS^2
\]
is holomorphic, thus it is rational. The unique preimage of $\infty$ by $F$ is itself, thus $F$ is a polynomial. It's degree must coincide with it's ramified covering degree, that is $2$.

The map $\phi\circ S^{-1}$ from $S(\bbD)$ to $\bbD$ is holomorphic (the Beltrami form $0$ is mapped to $\mu_0$ by $S^{-1}$ and then back to $0$ by $\phi$), maps $0$ to $0$, and conjugates $F$ to the rotation $R_\theta : \bbD \to \bbD$,
thus $F$ fixes $0$ with multiplier $e^{i 2\pi \theta}$. Thus $F(z) = \rho z+a z^2$ for some $a\in \bbC^*$. On one hand this map has only one critical point $z = -\rho/2a$, on the other hand, $S$ sends the ramification point $z=1$ of $\widetilde{B}$ to $z = -\rho/2$. Thus $a=1$:
\[
 F(z) = \rho z + z^2
\]
This proves that $\widetilde{B}$ is a model for $F$.


\section{The geometric construction of the pre-model}

The Blaschke product $B$ is usually introduced by a formula. If we want a pre-model for other families than the quadratic polynomials, it seems we need to find the corresponding formulas. We present here a geometric construction of $B$ that does not require a formula, and that easily generalizes. We carry out this generalization for some entire maps and also for some Lavaurs maps. In each case, the pre-model will be called $\beta$.

\subsection{For quadratic polynomials}

Let $D$ be the open ball in $\bbC$ of center $1$ and radius $1$.
\begin{center}
 \scalebox{0.7}{
 \begin{picture}(135,135)
  \put(0,0){%
   \scalebox{0.6}{%
    \includegraphics*[195,285][420,510]{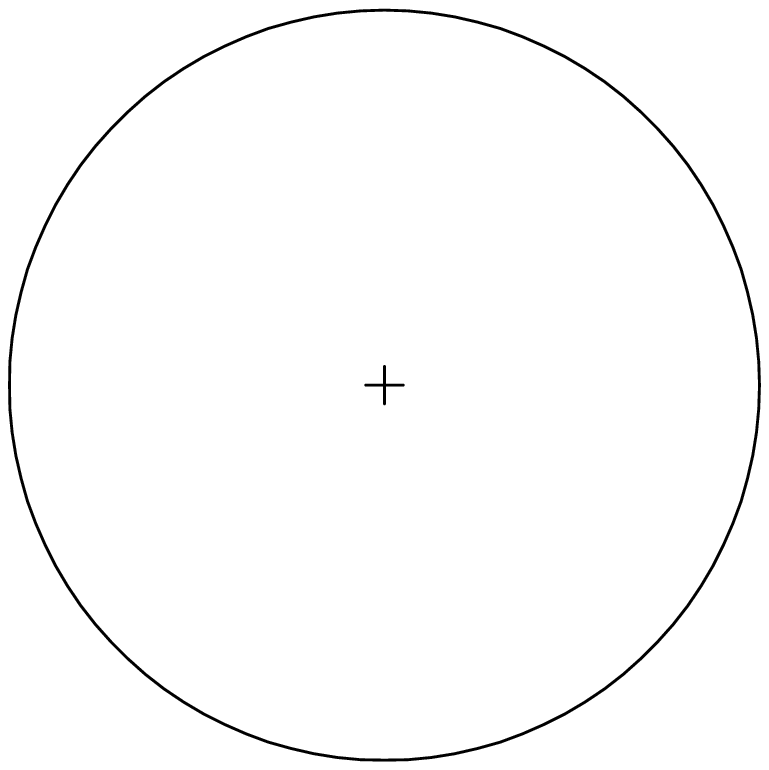}
   }
  }
  \put(-8,60){$0$}
  \put(75,57){$1$}   
  \put(90,90){$D$}   
 \end{picture}
 }
\end{center}
The preimage of $\partial D$ by the map $q : z \mapsto z^2$
is a $\infty$ shaped curve (a lemniscate).
\begin{center}
 \scalebox{0.7}{
 \begin{picture}(221,84)
  \put(0,0){%
   \scalebox{0.7}{%
    \includegraphics*[150,335][465,455]{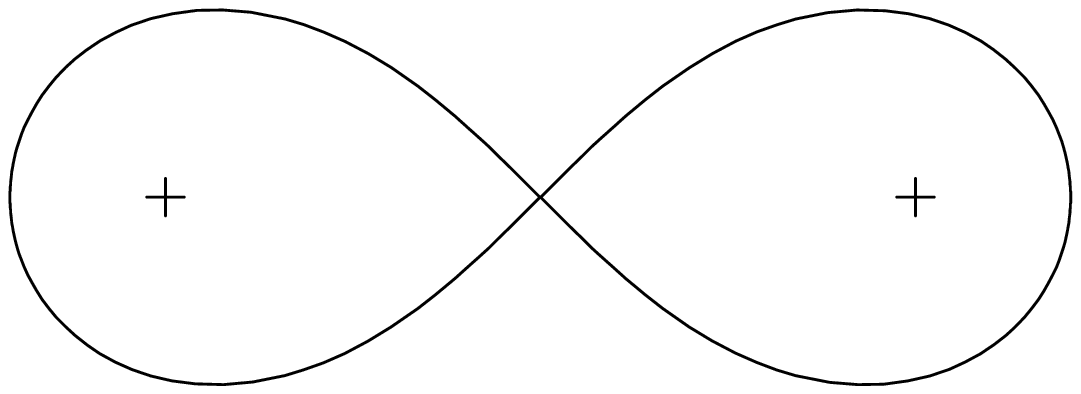}
   }
  }
  \put(30,30){$-1$}
  \put(107,30){$0$}
  \put(178,30){$1$}
 \end{picture}
 }
\end{center}
It cuts the complex plane into three
connected components. The unbounded one maps $2$-to-$1$ to $\bbC
\setminus \adh{D}$ by $q$, the other two map isomorphically to $D$.
Let us call $U$ the bounded one that contains $z=-1$, and let $W$ be
the complement of $\adh{U}$ in $\bbS^2$.
\begin{center}
 \scalebox{0.7}{
 \begin{picture}(221,84)
  \put(-2,0){%
   \scalebox{0.97}{%
    \includegraphics{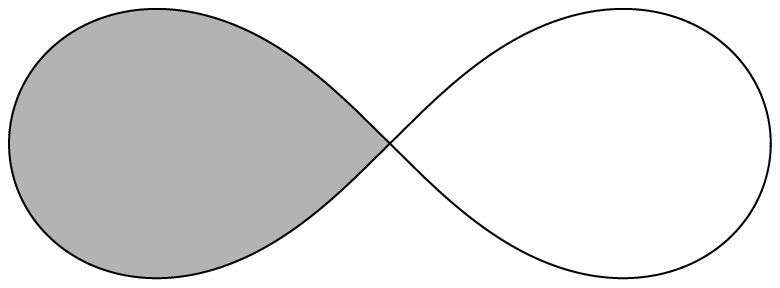}
   }
  }
  \put(30,35){$U$}
  \put(-30,35){$W$}
 \end{picture}
 }
\end{center}
Let $\phi$ be the unique
isomorphism from $W$ to $\bbS^2 \setminus \adh{\bbD}$ that fixes
$\infty$ and maps $0$ (the angular point of $U$) to $1$.
\begin{center}
 \scalebox{0.7}{
 \begin{picture}(221,100)
  \put(0,0){%
   \includegraphics{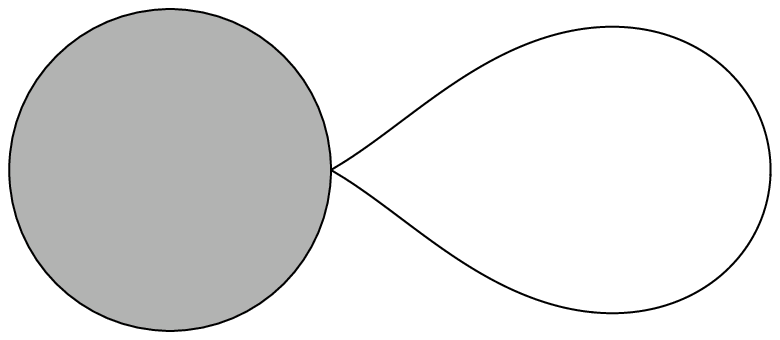}
  }
  \put(50,45){$\bbD$}
  \put(-35,45){\scalebox{1.5}{$\overset{\phi}{\longrightarrow}$}}
 \end{picture}
 }
\end{center}
Let $b = t \circ q \circ \phi^{-1} : \bbS^2 \setminus \adh{\bbD} \to \bbS^2$, where $t(z)=1-z$. Then $b(z)$
tends to $\partial \bbD$ when $z$ tends to $\partial \bbD$, thus $b$ admits a Schwarz
reflection along $\bbS^1$: that means there exists a
holomorphic function $B : \bbS^2 \to \bbS^2$ such that $B(z)=b(z)$
when $z\not\in \adh{\bbD}$, and such that $B(z) = s(b(s(z)))$ when
$z\in\bbD$, where $s(z) = 1/\bar{z}$. Moreover $B(\bbS^1) \subset
\bbS^1$, thus $B$ is a Blaschke fraction.
\begin{center}
 \begin{picture}(360,100)
  \put(0,0){%
   \includegraphics{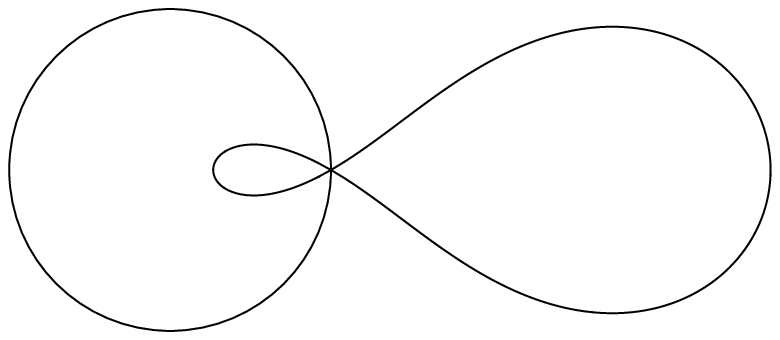}
  }
  \put(235,50){$\overset{B}{\longrightarrow}$}
  \put(258,0){%
   \scalebox{0.43}{%
    \includegraphics{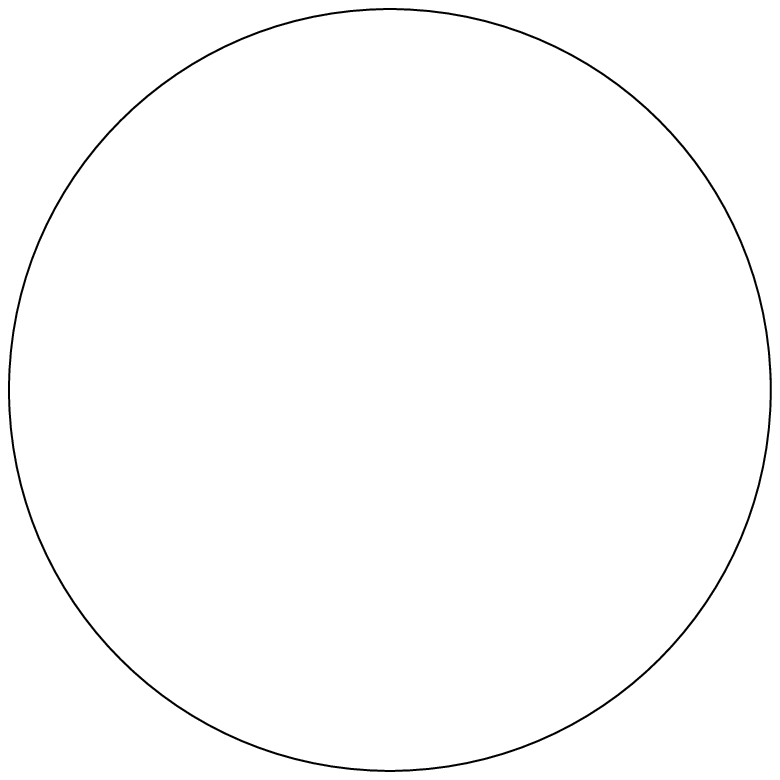}
   }
  }
 \end{picture}
\end{center}
Every point not in $\bbS^1
\cup \{0,\infty\}$ has $3$ preimages by $B$, thus $B$ has degree
$3$. It fixes $\infty$ with local degree $2$, thus it is of the form
\[B(z) = \lambda z^2 \frac{z-a}{1-\bar{a}z}\]
where $a \in \bbC^* \setminus \bbS^1$, and $\lambda\in \bbS^1$.
The point $z=1$ is the unique preimage of $1$, thus $z=1$ is a local
degree $3$ critical point of $B$. This implies $a = 3$. Finally, $1$
is mapped to $1$, thus $\lambda = 1$. We have proved that
\[B(z) = z^2 \frac{z-3}{1-3z}\]
The pre-model $\beta$ is just the function $B$.

\subsection{For some entire maps}\label{subsec_const_entire}

Let $\theta$ be a bounded type irrational and
$\rho=\exp(i 2\pi\theta)$. Let us give the class of entire maps we are able to study.

\varthm{Definition}
We assume $f : \bbC \to \bbC$ is
holomorphic, fixes $0$ with multiplier $\rho$, and has set of singular values contained in
$\{0,1\}$ and containing $\{1\}$.
We make one more assumption: that the connected component $U$ of
$f^{-1}(\bbD)$ containing $0$ is bounded.
\medbreak
This assumption is equivalent to the following: for any Jordan arc $\gamma$ from $0$ to $1$, the lift by $f$ of $\gamma$ starting from $0$ does not tend to $\infty$ (this implies it tends to some point, independent of the choice of $\gamma$, which we will call the \emph{main} critical point).

Note that we allow polynomials. In this case, the assumption is automatic.

These are very restrictive conditions. However, we will see that there
are uncountably many topologically inequivalent such maps in section~\ref{subsec_ex_entire}. 

The restriction $f : U \setminus f^{-1}(0) \to \bbD\setminus \{0\}$ is a covering. Since $U$ contains $0$ and $f$ has local degree $1$ at this point, $f : U \to \bbD$ is an isomorphism.
Moreover, $U$ is a Jordan disk (the bounded component of the complement of a Jordan curve), and $f : \adh{U} \to \adh{\bbD}$ is a homeomorphism.
Let $z$ be the preimage of $1$ by this homeomorphism.
Let $\phi$ be the conformal isomorphism between $\bbC\setminus \adh{U}$ and $\bbC\setminus \adh{\bbD}$ whose extension to boundaries map $w$ back to $1$. Consider the map $f \circ \phi^{-1} : \bbC\setminus\adh{\bbD} \to \bbC$. It tends to $\partial \bbD$ when the variable tends to $\partial \bbD$, thus it has a Schwarz reflection extension $\beta : \bbC\setminus\{0\} \to \bbC$.
The map $\beta$ is our pre-model.

\subsection{For some Lavaurs maps}\label{subsec_premodelLavaurs}

I have invented the geometric construction of pre-models in my thesis~\cite{C}, precisely because none were known for Lavaurs maps.

We will assume that the reader is familiar with the limit objects that arise in parabolic implosion. They were introduces, among others, by Douady, Hubbard, Lavaurs, Shishikura, \ldots We give a small summary in section~\ref{sec_obj}. I also give definitions and properties in~\cite{C} (in French). Some references are also~\cite{D},~\cite{Sh}.

We consider the quadratic polynomial
\[
 P=\upsilon z + z^2
\]
where $\upsilon$ is a root of unity. Let us write
\[
 \upsilon=\exp\big(i 2\pi \frac{p}{q}\big)
\]
with $q\in\bbN^*$, $p\in\bbZ$, and $p/q$ irreducible.
Let $K$ be the filled-in Julia set of $P$.

\begin{figure}[h]
 \begin{center}
 \begin{picture}(290,230)
  \put(0,30){
   \put(0,0){%
    \includegraphics*[10,30][300,220]{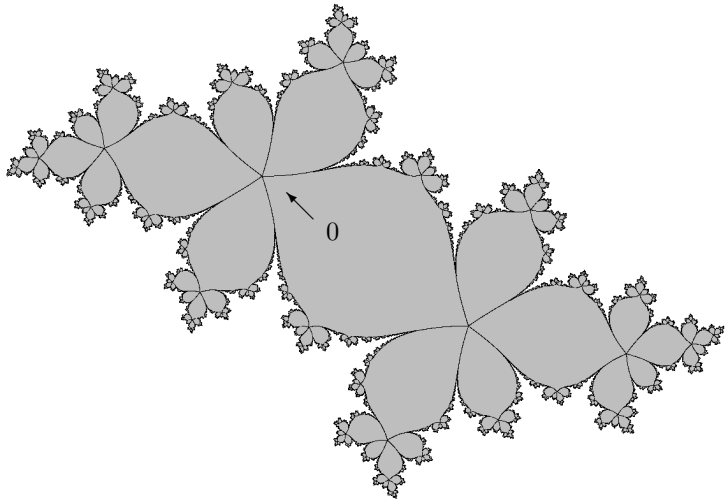}%
   }
   \put(130,100){$0$}
   \put(125,108){\vector(-1,1){10}}
  }
 \end{picture}
 \end{center}
 \caption{The filled-in Julia set of $P$ for $p/q =2/5$ (boundary in black, interior in gray).}
\end{figure}

Throughout this article, the notation $T_a$ refers to the translation by vector $a$ on either the real line $\bbR$, the complex plane $\bbC$, or their quotient by the subgroup $\bbZ$.

The repelling axes of the parabolic point $z=0$ are labeled by $j$, which can take $q$ values.
The attracting ones are labeled by $i$, which can also take $q$ values.
The extended repelling Fatou parameterization associated to the axis $j$ is noted
\[
 \psi_{+,j} : \bbC_j \to \bbC
\]
where $\bbC_j$ is just a copy of $\bbC$ that we labeled $j$. It is a semi-conjugacy from $T_1$ to $P^q$: $\psi_{+,j} \circ T_1 = P^q \circ \psi_{+,j}$.
The \emph{horn maps} are maps
\[
 h_{j,\sigma} : \psi_{+,j}^{-1}\ins{K} \to \bbC_j
\]
Their domain of definition does not depend on $\sigma$. It is invariant by $T_1$, and $h_{j,\sigma}$ commutes with the translation $T_1$. Thus there is a quotient map from a subset of the cylinder $\bbC/\bbZ$ to the cylinder, that we will note
\[
 \overline{h}_{j,\sigma}
\]
This quotient is defined in neighborhoods of the the two ends of the cylinder, and has a holomorphic extension to these ends. This extension will be noted
\[
 \widehat{h}_{j,\sigma}
\]  
It fixes both ends, with non zero multiplier.

For a fixed $\sigma$, the $q$ maps $\widehat{h}_{j,\sigma}$ with different $j$ are all conjugated by a translation. Thus now, we will omit the index $j$, because to study their dynamics, it is enough to work with any one of them.

Let $\theta$ be a bounded type irrational number. From the definition of the horn maps, it follows that $h_{\sigma} = T_{\sigma} \circ h_{0}$. Thus there is a value of $\sigma$, unique modulo $\bbZ$, such that the upper end has multiplier $\exp (i 2\pi \theta)$. The map $\widehat{h}_{\sigma}$ has thus a Siegel disk there, and we are looking for a model.

The construction of the pre-model does not depend on $\sigma$. So, for the rest of this section, we choose any value of $\sigma$, and write $h$ in place of $h_\sigma$. The map $\overline{h}$ has only one critical value $v$. It is an infinite degree ramified cover of the cylinder. Back to the complex plane, the preimage by $h$ of the horizontal line $L$ through $v$ is called the \emph{chessboard graph}. The connected components of the complement in $\operatorname{Def}(h)$ of this preimage, are called the \emph{chessboard boxes}. Each is mapped isomorphically by $h$ to either the upper or the lower half plane delimited by the line $L$.

There is one chessboard box that is a neighborhood of the upper end of the cylinder. Let us call it $U$. Then $U$ is a Jordan disk in the cylinder completed by it's ends, and $\overline{h}$ is a homeomorphism from the closure of $U$ to the upper half cylinder delimited by $L$. There is thus exactly one point $w$ on the boundary of $U$ that maps to $v$. Let $W$ be the complement in the cylinder of the closure of $U$. Let $\phi$ be the unique analytic isomorphism between $W$ and the lower half cylinder delimited by the real line $\bbR$, that maps $w$ to $0$.

\newpage

\begin{figure}[h]
  \begin{center}
  \rotatebox{270}{%
   \includegraphics{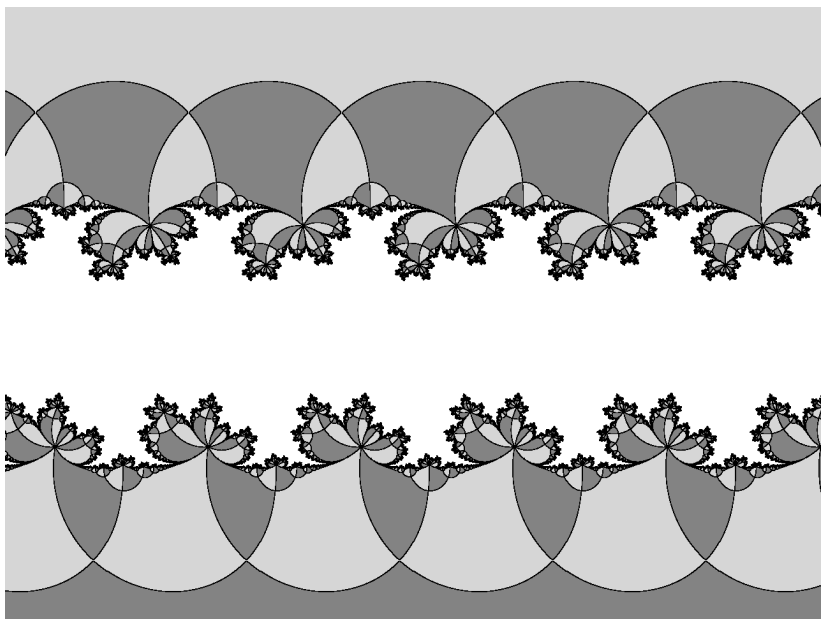}%
  }
  \end{center}
  \caption{Parabolic chessboard of $h$ for $p/q=2/5$; the boxes that map to the upper half plane are in light gray, those mapping to the lower half plane in dark gray. The white set is the preimage by $\psi_+$ of $\bbC\setminus K$. It is open and it's closure is the complement of $\operatorname{Def}(h)$.}
\end{figure}

\begin{figure}[h]
 \begin{center}
 \begin{picture}(264,176)
  \put(0,10){%
   \includegraphics{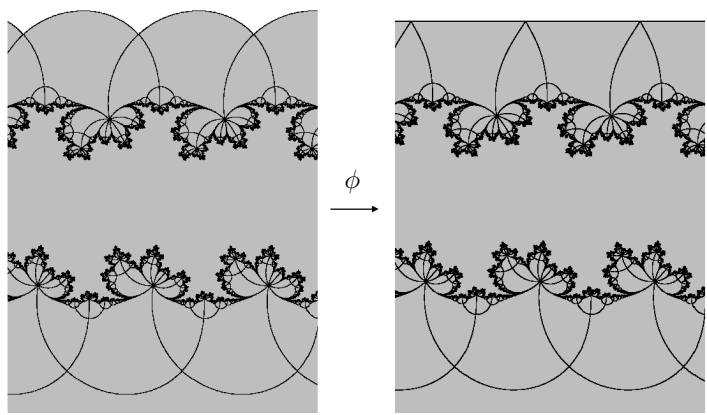}%
  }
  \put(128,95){$\phi$}
 \end{picture}
 \end{center}
 \caption{The map $\phi$ is an isomorphism of the gray parts.}
\end{figure}

\newpage

The map $T_{-v} \circ \overline{h} \circ \phi^{-1}$ has a Schwarz reflection $\beta_0$, because points tending to the real line have image tending to the real line. The map $\beta_0$ goes from a $T_1$-invariant subset of $\bbC$ to $\bbC$, and commutes with $T_1$. It's quotient $\beta$ is our pre-model: it goes from a subset of the cylinder $\bbC/\bbZ$ to $\bbC/\bbZ$.
The map $\beta_0$ sends $0$ to $0$, it's set of real critical points is $\bbZ$, and it's set of \emph{all} critical values is $\bbZ$. Thus $\beta$ has only one critical value.

\bigskip

\begin{center}
 \begin{picture}(351,327) 
  \put(0,10){%
   \includegraphics{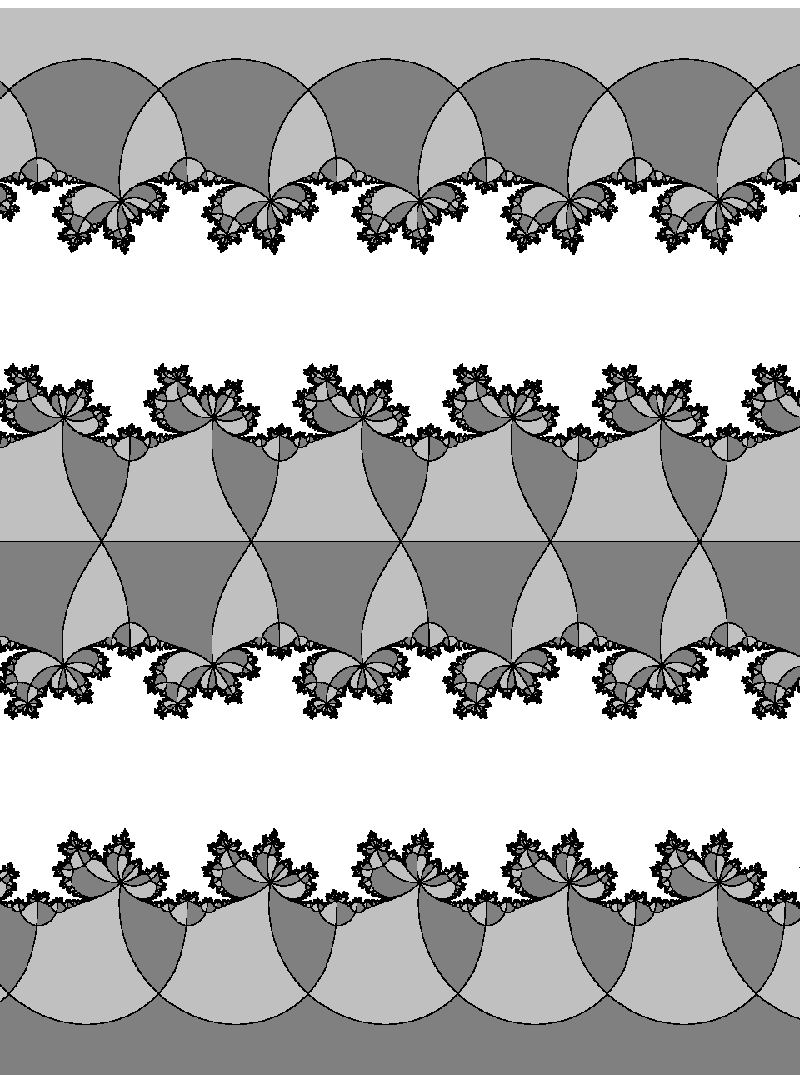}%
   }
  \put(270,2){\scalebox{0.35}{%
   \includegraphics{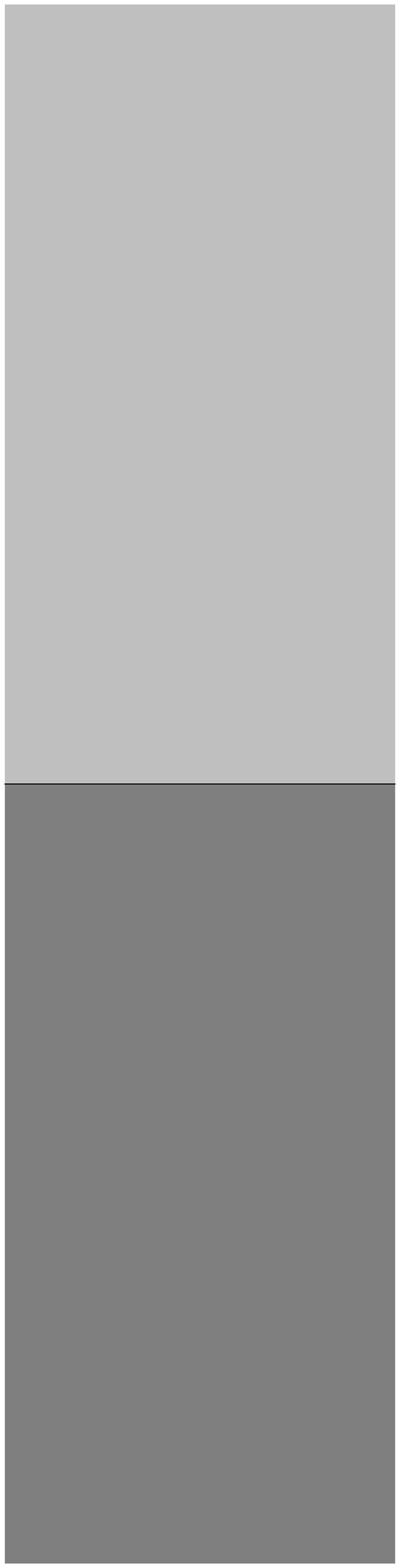}%
  }}
  \put(239,161){\scalebox{1.5}{$\overset{\beta_0}{\longrightarrow}$}}
 \end{picture}
\end{center}


\section{From the pre-model to the model}\label{sec_premodeltomodel}

This section describes a standard argument, due, among others, to Ghys, Douady, Shishikura, Herman and \'Swi\c{a}tek.
In this article, quasiconformal, quasiregular and quasisymmetric maps will be always understood to be orientation preserving. We assume the reader is familiar with these notions and with the straightening theorem (also called quasiconformal Riemann mapping theorem) of Beltrami forms (also called almost-complex structures, or ellipsis fields).
For an irrational $\theta$ with continued fraction expansion $\theta=a_0+1/(a_1+1/(a_2+\ldots))$, let us call $\sup_{i>0}a_i$ the \emph{type} of $\theta$.

\begin{theo}[Herman, \'Swi\c{a}tek]\label{thm_HeSw}
  Suppose $f: \bbT \to \bbT$, where $\bbT=\bbR/\bbZ$, is an orientation preserving homeomorphism such that
  \begin{enumerate}
   \item $f$ has rotation number $\theta$, irrational of bounded type, say bounded by $M\in\bbN^*$
   \item $f$ is $C^1$ on $\bbT$
   \item $f$ has a unique critical point, which is $0$
   \item there is a neighborhood $W$ of $0$ such that $f$ is $C^3$ on $W$, and has negative schwarzian derivative on $W\setminus\{0\}$
   \item $\exists l\in\bbN^*$ and real numbers $A,A'>0$, such that $\forall x\in W$, $A|x|^{2l} \leq f'(x) \leq A'|x|^{2l}$
   \item the variation of $\log f'$ on $\bbT\setminus W$ is bounded by $M'>0$
  \end{enumerate}
  Then, there is some $c>1$ such that $f$ is $c$-quasisymmetrically conjugate to the translation by $\theta$ on $\bbT$.
\end{theo}
I have borrowed this statement from \cite{Pe2}. On the other hand, if $\theta$ does not have bounded type,
then $f$ is never quasisymmetrically conjugate to the rotation $T_\theta$.

\

In the preceding section, we have constructed a pre-model map $\beta$, whose restriction to $\bbS^1$ (for entire maps) or $\bbT$ (for Lavaurs maps), is a homeomorphism. The quotient $\bbT$ is identified with the unit circle $\bbS^1$ in $\bbC$ by means of the map $t \mapsto \exp(i 2\pi t)$. The cylinder $\bbC/\bbZ$ is identified with the subset $\bbC^*$ of the Riemann sphere $\bbS^2$, by means of the same formula. Let $\theta$ be an irrational. There is a unique real number $\tau$ such that $R_\tau \circ \beta$ has rotation number $\theta$ on $\bbS^1$ (see \cite{dMvS}). Suppose now $\theta$ has bounded type. The map $R_\tau \circ \beta$ satisfies the conditions of the Herman-\'Swi\c{a}tek theorem. Let $\eta$ be the quasisymmetric conjugacy given by this theorem.

\begin{theo}[Ahlfors, Beurling]\label{thm_AhBe}
  Every quasisymmetric map $f : \bbS^1 \to \bbS^1$ extends to a
  quasiconformal map $\bbS^2 \to \bbS^2$. 
\end{theo}

Douady and Earle enhanced this theorem in~\cite{DE}.

\

Let $\widetilde{\eta}$ be such an extension of $\eta$ to $\bbD$. We will moreover require that $\widetilde{\eta}$ fix $0$, which is possible by post-composing any quasiconformal extension with a $C^\infty$ map equal to the identity in a neighborhood of the unit circle. Let $\mu_0$ be the pull-back of the null Beltrami form by $\widetilde{\eta}$: $\mu_0$ is defined on $\bbD$. Let $\widetilde{\beta} = R_\tau \circ \beta$ outside of $\bbD$, and $\widetilde{\beta} = \widetilde{\eta}^{-1} \circ R_\theta \circ \widetilde{\eta}$.
Thus it is possible to define globally a $\widetilde{\beta}$ invariant Beltrami form $\mu$, by means of the dynamics. The well known procedure goes as follows: if a point $z$ never falls in $\bbD$, let $\mu(z)$=0. Otherwise, let $z_n$ the first iterate that falls in $\bbD$. Then define $\mu(z)$ to be the pull-back of $\mu_0(z_n)$ by $\widetilde{\beta}^n$. 

The map $\widetilde{\beta}$ we just defined is called the model.


\section{Rigidity}

The notion of rigidity involved here is different from the usual notion of quasiconformal rigidity. In the classical rigidity theory one studies quasiconformal conjugacy classes of holomorphic maps. Here we focus on equivalence classes (which are simpler).

\subsection{Entire maps}

A subset of $\bbC$ will be called non accumulating if it has no accumulation points \emph{in $\bbC$}. This is equivalent to be closed (rel.\ $\bbC$) and discrete.
The two following statements are not hard to prove. However, I added the proofs. 

\begin{prop}[Rigidity]\label{prop_rigidity}
 Suppose $f$ and $g$ are two entire maps with a non accumulating set of singular values, and $\phi$ and $\psi$ two orientation preserving homeomorphisms of $\bbC$, such that $\psi$ is isotopic to identity rel.\ a non accumulating set $V$ containing all the singular values of $f$ and such that the following diagram commutes
 \[
  \xymatrix{
   \bbC \ar[r]^\phi \ar[d]_f & \bbC \ar[d]^g \\
   \bbC \ar[r]_\psi & \bbC \\
  }
 \]
 Then there exists a complex-affine function $a$ (i.e.\ of the form $a(z)=b z + c$) such that the following diagram commutes 
 \[
  \xymatrix{
   \bbC \ar[r]^a \ar[d]_f & \bbC \ar[d]^g \\
   \bbC \ar[r]_{\operatorname{id}} & \bbC \\
  }
 \]
and such that for all $z\in f^{-1}(V)$, $a(z) = \phi(z)$.
\end{prop}
\varthm{Proof}
 Let $u : [0,1]\times \bbC$ be the isotopy and let us note $u_t(z)=u(t,z)$:
 \begin{itemize}
  \item $u_0=\psi$,
  \item $u_1=\operatorname{id}_\bbC$,
  \item $\forall t\in [0,1]$ and $z\in V$, $u_t(z)=z$.
 \end{itemize}
 We would like to complete the following diagram for all $t$:
 \[
  \xymatrix{
   \bbC \ar@{.>}[r]^{v_t} \ar[d]_f & \bbC \ar[d]^g \\
   \bbC \ar[r]_{u_t} & \bbC \\
  }
 \]
 Let $Y = \bbC\setminus f^{-1}(V)$ and $Z = \bbC\setminus g^{-1}(V)$. We have the following commutative diagram:
 \[
  \xymatrix{
   Y \ar[r]^{\phi} \ar[d]_{f} & Z \ar[d]^{g} \\
   \bbC \setminus V \ar[r]_{\psi} & \bbC \setminus V \\
  }
 \]
Note that g is a covering from $Z$ to $\bbC \setminus V$. Let us consider the function $\zeta(t,z) = \zeta_t(z) = u_t \circ f (z)$. Since $u_t(z)=z$ whenever $z$ belongs to $V$, this implies $\zeta^{-1}( \bbC \setminus V) = [0,1] \times Y$. We want to complete the following diagram:
 \[
  \xymatrix{
   [0,1]\times Y \ar[dr]^{\zeta} \ar@{.>}[r]^v & Z \ar[d]^{g} \\
   & \bbC \setminus V \\
  }
 \]
Let us choose a base point $z_0\in Y$ and map it by $\phi$ to a basepoint of $Z$. Let $(0,z_0)$ be the basepoint of $[0,1] \times Y$. With the further requirement that $v(0,z_0) = \phi(z_0)$, the lifting theorem states that $v$ exists if and only if
\[
\zeta_*(\pi_1 ([0,1] \times Y)) \subset g_*(\pi_1 Z)
\]
in which case $v$ is unique. By retracting $[0,1] \times Y$ to $\{0\} \times Y$, we have $\zeta_*(\pi_1 ([0,1] \times Y)) = (\psi\circ f)_*(\pi_1 Y) = (g \circ \phi)_*(\pi_1 Y) \subset g_*(\pi_1 Z)$.
Thus $v$ exists, and coincides with $\phi$ for $t=0$. 

 Each map $v_t$ is continuous. We claim that the functions $v_t$ extend to homeomorphisms of the plane by sending every point in $z \in f^{-1}(V)$ to $\phi(z)$, which is independent of $t$. Let us first prove injectivity on $Y$. If $v_t$ takes the same value on two points $z$ and $z'$ in $Y$, then the commutative diagram implies that $u_t\circ f$ takes the same value, and since $u_t$ is injective, $f(z) = f(z')$. Thus for all $s$, $g \circ v_{s} = u_{s} \circ f$ takes the same value. Now since $g$ is a covering from $Z$ to $\bbC\setminus V$, and $v_t(z)=v_t(z')$, this implies that $v_{s}(z)=v_{s}(z')$ holds for all $s$. In particular it holds for $s=0$. But $v_0 = \phi$ is injective, thus $z=z'$, q.e.d. Now $v_t$ is the lift of the covering $\zeta_t = u_t \circ f$, thus it is a covering ($\bbC \setminus V$ is connected and locally path connected). Injective coverings (from a non empty space to a connected one) are homeomorphisms, and so $v_t$ is a homeomorphism from $Y$ to $Z$.
Let us extend $v_t$ on all of $\bbC$ by setting $v_t(z)=\phi(z)$ when $z \in \bbC \setminus Y$. These extensions are bijections of $\bbC$ that make the diagram commute. An homotopy argument with a small injective loop around an element of $\bbC\setminus Z$ and winding number considerations shows that these extensions are homeomorphisms.

Let $a(z) = v_1(z)$ for $z \in \bbC$. Then $a$ is continuous and $g\circ a = f$. Thus $a$ is holomorphic. We saw that it is also a homeomorphism of $\bbC$. Thus this is affine.
\qed\medbreak

\begin{cor}\label{cor_topeqaffeq}
 If $f$ and $g$ are two entire maps with at most two singular values and if $f$ and $g$ are positively topologically equivalent, then they are affine equivalent. Moreover, let $V$ be a pair of distinct points containing the critical values of $f$ and let us write respectively the topological equivalence and the affine one with
\[
 \xymatrix{
   \bbC \ar[r]^u \ar[d]_f & \bbC \ar[d]^g \\
   \bbC \ar[r]_v & \bbC \\
 }
 \qquad\qquad
 \xymatrix{
   \bbC \ar[r]^a \ar[d]_f & \bbC \ar[d]^g \\
   \bbC \ar[r]_b & \bbC \\
 }
\]
Then $a$ and $b$ can be chosen so that $v$ and $b$ take the same value on the set $V$, and such that $u$ and $a$ take the same values on $f^{-1}(V)$. 
\end{cor}
\varthm{Proof}
 Let us note $V = \{z_0,z_1\}$. Let $A$ be the unique affine map mapping $z_0$ to $0$ and $z_1$ to $1$. Let $w_0=v(z_0)$ and $w_1=v(z_1)$. Let $B$ be the unique affine map mapping $w_0$ to $0$ and $w_1$ to $1$. Let $F=A \circ f$ and $G=B\circ g$. Then $F $ and $G$ satisfy the hypotheses of proposition~\ref{prop_rigidity}, with $\phi=u$, $\psi = B \circ v \circ A^{-1}$ and $V'=A(V)$: indeed, $\psi$ fixes $0$ and $1$, and the only non trivial hypothesis to check is that $\psi$ is isotopic to identity rel.\ $\{0,1\}$. A topological theorem states that every orientation preserving homeomorphism of $\bbC$ is isotopic to identity. Let $t \mapsto \psi_t$ be such an isotopy, with $t\in [0,1]$, $\psi_0 = \psi$ and $\psi_1 = \operatorname{id}$. Let $C_t$ be the unique affine homeomorphism mapping $\psi_t(0)$ to $0$ and $\psi_t(1)$ to $1$. Then $C_t$ varies continuously with $t$, and the map $t\mapsto C_t \circ \psi_t$ is an isotopy rel.\ $\{0,1\}$ from $\psi$ to $\operatorname{id}$.  So we can apply proposition~\ref{prop_rigidity}, and this tells us that there exists some affine map $a$ such that $F = G \circ a$, i.e.\ $A \circ f = B \circ g \circ a$, i.e.\ $b \circ f = g \circ a$ where $b = (B^{-1} \circ A)$. Proposition~\ref{prop_rigidity} also states that $a$ and $u$ have the same values on the set $F^{-1}(V') = f^{-1}(V)$.
 It is easy to check that $b$ and $v$ take the same values on $z_0$ and $z_1$.
\qed\medbreak

\subsection{Lavaurs maps}

The difference here is that these maps are not rigid~! Even in the cylinder. Take any
homeomorphism $u : \bbC \to \bbC$ which is conformal on the domain of definition of $h$, commutes with $T_1$, but changes the $\bbC$-affine shape of $\operatorname{Def}(h)$. Then the map $h \circ u$ is equivalent (not conjugated) to $h$, but not affine equivalent.

To recover rigidity we need further assumptions on the maps in the equivalence, and we will make use of the quasiconformal conjugacy rigidity of the polynomial $P(z)=\upsilon z+z^2$.

\begin{prop}[folk.]\label{prop_folk}
 Let $\upsilon$ be a root of unity. Any holomorphic map $f : \bbC \to \bbC$ quasiconformally conjugated to $P=\upsilon z + z^2$ is affine conjugated to $P$. Moreover, if $f=P$, then the conjugacy fixes any point in the grand orbit of the critical point, and any point in the Julia set $J$, and it leaves invariant each component of the complement of $J$.
\end{prop}
In fact, it is known that the quasiconformal conjugacy rigid quadratic polynomials are exactly those who belong to the boundary of the Mandelbrot set, and the above proposition holds for them.

Consider now a Beltrami form $\mu$ that is the pull-back by $\psi_+$ of some Beltrami form $\mu'$ on $\bbC$. Then $\mu$ is $T_1$-invariant if and only if $\mu'$ is $P^q$-invariant.
\begin{prop}\label{prop_mucomp}
 Let $P$ be as in the preceding proposition.
 Let $\mu$ be a Beltrami form defined on $\bbC$ such that:
 \begin{itemize}
  \item $\mu$ is the pull-back by $\psi_+$ of some $P$-invariant Beltrami form $\mu'$ on $\bbC$.
  \item $\|\mu\|_\infty<1$
 \end{itemize}
 Let $J$ be the Julia set of $P$.
 Let $J'=J \cup$ the grand orbit of the critical point of $P$.
 Let $w_0$ be any point in $\psi_+^{-1}(J')$.
 Let $S$ be the unique straightening of $\mu$ which commutes with $T_1$ and fixes $w_0$.
 Then $S$ fixes any point in $\psi_+^{-1}(J')$, and leaves invariant any component of the complement of $\psi_+^{-1}(J)$.
\end{prop}
Note that the hypotheses imply that $\mu = \lambda \mu_0$ on the complement of $\operatorname{Def}(h)$, where $\lambda\in\bbC$, $|\lambda|<1$ and $\mu_0$ is the pull-back by $\psi_+$ of the line field parallel to external rays of the Julia set. There is much more freedom in $\operatorname{Def}(h)$ (the vector space of all such $\mu$'s is infinite dimensional). Note also that we not only require $\mu'$ to be $P^q$-invariant, but also $P$-invariant.
\varthm{Proof}
 We just give an idea of the proof.
 According to proposition~\ref{prop_folk}, there is a straightening $R$ of $\mu'$ that conjugates $P$ to itself. We claim that there is a map $S : \bbC \to \bbC$ such that $\psi_+ \circ S = R \circ \psi_+$. To prove this, one considers the set $\mcB$ of backward orbits of $P^q$ (sequences $(z_n)_{n\in\bbN}$ with $P^q(z_{n+1}) = z_n$) that tend to $0$ along the repelling axis of index $j$ (the index that we fixed in section~\ref{sec_obj}). The map that sends $w$ to the sequence $\psi_+(w-n)$ is a bijection from $\bbC$ to $\mcB$. Moreover $R$ acts on $\mcB$, from which the claim follows. The same can be done for $R^{-1}$, which provides an inverse for $S$. One then proves that $S$ and $S^{-1}$ are continuous, thus $S$ is a homeomorphism. (For this, let $V$ be a small disk centered on $0$, on which $P^q$ has an inverse branch. Every element of $(z_n)\in\mcB$ has an index $N$ for which $\forall n\geq N$, $z_n \in V$. Let $w_k$ be a sequence in $\bbC$, and $z_n^k = \psi_+(w_k-n)$. Let $w\in\bbC$ and $z_n=\psi_+(w-n)$. Then it is enough to notice that $w_k$ tends to $w$ if and only if there exists a common index $N$ for the sequences $z_n$, $z_n^k$ and $z_N^k \tend z_N$ as $k \to +\infty$.) The relation $\psi_+ \circ S = R \circ \psi_+$ implies that $S$ is quasiconformal and maps $\psi_+^*(\mu')=\mu$ to $\psi_+^*(0)=0$: in other words it straightens $\mu$. Finally, $R$ leaves invariant any element $(z_n) \in \mcB$ such that $z_0\in J'$, and the last statement of the proposition follows.
\qed\medbreak

The proof yielded the following
\begin{prop}\label{comp_R}
 Let $R$ be the straightening of $\mu'$ that conjugates $P$ to itself:
 \[
  \psi_+ \circ S = R \circ \psi_+
 \]
\end{prop}

We will only make use of one implication of the next proposition, but we think it is worth to mention the equivalence.

\begin{prop}\label{prop_pullpsih}
 For a Beltrami form $\mu$ on $\bbC$ that is equal to $0$ on $\bbC\setminus \operatorname{Def}(h)$, the condition that $\mu$ be a pull-back by $\psi_+$ of a $P$-invariant Beltrami form is equivalent to the restriction of $\mu$ to $\operatorname{Def}(h)$ being a pull-back by $h$ of a $T_1$-invariant Beltrami form.
\end{prop}
\varthm{Proof}
 The map $h$ has the form $h= T_{\operatorname{sg}} \circ \phi_\div \circ \psi_+$ (see section~\ref{sec_obj} for the definition of $\phi_\div$), so the proposition amounts to: the Beltrami forms $\mu'$ on $\insb{K}$ that are $P$-invariant are exactly the pull-backs by $\phi_\div$ of $T_1$-invariant forms. This claim follows from the remark that two points $z$ and $z'$ belong to the same grand-orbit by $P$ if and only if $\phi_\div(z) = \phi_\div(z) \bmod \bbZ$.
\qed\medbreak

\begin{prop}\label{prop_l}
 Suppose $l : \operatorname{Def}(l) \to \bbC$ is holomorphic with $\operatorname{Def}(l) \subset \bbC$.
 Suppose $u$ and $v$ are quasiconformal homeomorphisms of $\bbC$, that both commute with $T_1$, that $u$ fixes at least one point $w_0$ in $\psi_+^{-1}(J')$ (where $J'=J \cup$ the grand orbit of the critical point of $P$), that $\overline{\partial} u =0$ on $\bbC\setminus \operatorname{Def}(h)$, that $u(\operatorname{Def}(h)) = \operatorname{Def}(l)$, and that the following diagram commutes
 \[
  \xymatrix@!{
   \operatorname{Def}(h) \ar[r]^u \ar[d]_h & \operatorname{Def}(l) \ar[d]^l \\
   \bbC \ar[r]_v & \bbC \\
  }
 \]
Then $\operatorname{Def}(l) = \operatorname{Def}(h)$ and there exists $\kappa \in \bbC$ such that
 \[
  l = T_\kappa \circ h
 \]
\end{prop}
\varthm{Proof}
First, the hypotheses imply that there is a $\kappa\in\bbC$ such that $v$ coincides with $T_\kappa$ on the set of critical values of $h$.
Let $\mu' = u^*(0)$, and $\mu = v^*(0)$. The Beltrami form $\mu'$ is null on $\bbC\setminus \operatorname{Def}(h)$ and, on $\operatorname{Def}(h)$, it is the pull-back by $h$ of $\mu$. Thus according to proposition~\ref{prop_pullpsih}, proposition~\ref{prop_mucomp} applies to $\mu'$. In particular, $u$ preserves each component $U$ of $\operatorname{Def}(h)$.
Now, for $t\in [0,1]$, let $v_t$ be the straightening of $t.\mu$ that coincides with $v$ on the set of critical values of $h$: $v_1 = v$ and $v_0 = T_\kappa$. Let $\pi : \bbC \to \bbC/\bbZ$ be the quotient map. Since $\pi \circ h$ is a ramified cover (of infinite degree), with only one critical value, an isotopy lifting argument, like in the preceding section, shows that the isotopy $v_t$ can be pulled-back through the diagram into an isotopy $u_t$ on $\operatorname{Def}(h)$ such that $l\circ u_t = v_t \circ h$, $u_t$ commutes with $T_1$, and $u_1 = u\big|_{\operatorname{Def}(h)}$. The map $u_0$ is an analytic isomorphism of each component $U$ of $\operatorname{Def}(h)$. Since $u_t$ maps a critical point $z$ of $h$ to a critical point of $l$, $u_t(z)$ can not vary with $t$. Thus $u_0(z)=u_1(z)=u(z)=z$ for all critical points $z$ of $h$. Now, each component of $\operatorname{Def}(h)$ contains at least $2$ critical points of $h$, in fact infinitely many, and has hyperbolic universal covering (a component omits way more than two points in $\bbC$), thus $u_0$ is the identity on $\operatorname{Def}(h)$. This proves the proposition.
\qed\medbreak


\section{The straightening}

Here, we prove that straightening the Beltrami form provides the expected map.

\subsection{Entire maps}

Let us consider a map $f$ as in section~\ref{subsec_const_entire}. Let $\beta$ be the associated pre-model. Let $\theta$ be any bounded type irrational and $\widetilde{\beta}$ the model defined in section~\ref{sec_premodeltomodel}.
A conformal isomorphism $\phi : \bbC\setminus \adh{U} \to \bbC\setminus \adh{\bbD}$ was defined, where $U$ is the connected component of $f^{-1}(\bbD)$ that contains $0$. There was also a homeomorphism $\widetilde{\eta}$ of $\adh{\bbD}$, fixing $0$ and quasiconformal on $\bbD$. The map $\widetilde{\beta}$ was defined by
\[
 \widetilde{\beta} = R_\tau \circ f \circ \phi^{-1}
\]
on $\bbC\setminus \adh{\bbD}$ and
\[
 \widetilde{\beta} = \widetilde{\eta}^{-1} \circ R_\theta \circ \widetilde{\eta}
\]
on $\adh{\bbD}$.
Let us extend $\phi$ into a (quasiconformal) homeomorphism $\widetilde{\phi}$ of $\bbC$, by means of the formula
\[
  \widetilde{\phi} = \widetilde{\eta}^{-1} \circ R_{-\theta} \circ \widetilde{\eta} \circ R_\tau \circ f
\]
on $\adh{U}$, with the notations of section~\ref{sec_premodeltomodel}.
This extension has been chosen so that the relation
\[
 \widetilde{\beta} \circ \widetilde{\phi} = R_\tau \circ f
\]
holds on all of $\bbC$. Thus $\widetilde{\beta}$ and $f$ are topologically equivalent.
Let $\mu$ the $\widetilde{\beta}$-invariant Beltrami form defined in section~\ref{sec_premodeltomodel}, $S$ the unique straightening (sending $\mu$ to $0$) that maps $0$ to $0$ and that maps the unique non-zero singular value of $\widetilde{\beta}$, namely $R_t(1)$, to $1$,
and $g$ the conjugated function
\[
 g = S \circ \widetilde{\beta} \circ S^{-1}
\]
It is holomorphic, entire, and topologically equivalent to $f$ with the following equivalence
\[
 g \circ (S \circ \widetilde{\phi}) = (S \circ R_t) \circ f
\]
According to corollary~\ref{cor_topeqaffeq}: $g\circ a = b\circ f$ with two affine isomorphisms $a$ and $b$, and with $b$ taking on $0$ and $1$ the same value as $S\circ R_\tau$, i.e.\ $0$ and $1$, thus $b=\operatorname{id}$, and with $a$ taking the same value as $S \circ \widetilde{\phi}$ at every point of $f^{-1}(\{0,1\})$.
Thus $a(0)=S \circ \widetilde{\phi}(0)=0$, which implies $a$ is of the form $z \mapsto \lambda z$ for some $\lambda \in \bbC^*$. Now the equivalence amounts to $g(\lambda z) = f(z)$. Thus $\lambda g'(0) = f'(0)$.
But near $0$, $g$ is conjugated to $R_\theta$ by the map $\widetilde{\eta} \circ S^{-1}$. This conjugacy is conformal because the pull-back of the null Beltrami form by $\widetilde{\eta}$ is $\mu$ and $S$ straightens $\mu$. Thus $g'(0) = \exp(i 2\pi \theta)$. Thus $g'(0) = f'(0)$, which implies $\lambda = 1$. We have proved that
\[
 g = f
\]
In other words: the model $\widetilde{\beta}$ is quasiconformally conjugated to $f$.

\begin{cor}
  For all entire map satisfying the conditions of section~\ref{subsec_const_entire}, there is a Siegel disk containing $0$, whose boundary is a quasicircle going through a (unique) critical point.
\end{cor}

This critical point must be the one I called the main critical point in section~\ref{subsec_const_entire}.

\subsection{Lavaurs maps}

Let us sum up what we have obtained up to now, concerning Lavaurs maps. We had a root of unity $\upsilon$, the polynomial $P=\upsilon z +z^2$, we chose one of the repelling axes $j$, and denoted $\psi_+$ the repelling Fatou coordinate from $\bbC_j$ to $\bbC$. We constructed a pre-model map $\beta$. We were then given a bounded type irrational $\theta$ and we associated to it a model map $\widetilde{\beta}$. These objects (sets, maps, Beltrami form) live in the cylinder $\bbC/\bbZ$, and we now want to consider their lifts to the universal cover $\bbC$. But for convenience, \emph{we will use the same notations}. We constructed a $\widetilde{\beta}$-invariant and $T_1$-invariant Beltrami form $\mu$. Let $w_0$ be any point in $\psi_+^{-1}(J)$. Straightening $\mu$ conjugates $\widetilde{\beta}$ to a holomorphic map $l$:
\[
 l = S \circ \widetilde{\beta} \circ S^{-1}
\]
where $S$ is the unique straightening of $\mu$ that commutes with $T_1$ and sends $\phi(w_0)$ to $w_0$ (the map $\phi$ is defined is section~\ref{subsec_premodelLavaurs}; the condition on $w_0$ is here to ensure that the domain of definition of $l$ coincides with that of $h$, as we will prove in a few paragraphs). The following theorem states that we have obtained what we were looking for:

\begin{theo}\label{thm_lish}
 \[ \exists \sigma \in \bbC\ \ l = h_\sigma \]
 and
 \[(\widehat{h}_\sigma)'(+i \infty) = \exp(i 2\pi \theta)\]
\end{theo}

So, $\widetilde{\beta}$ is a quasiconformal model for the horn map with a Siegel disk with rotation number $\theta$ at the upper end of the cylinder.

\

\noindent \emph{\underline{Proof of theorem~\ref{thm_lish}}}:

\

\noindent We want to apply propositions~\ref{prop_mucomp} and~\ref{prop_l}.
The map $\widetilde{\beta}$ and the Beltrami form $\mu$ are not well suited, because their domain of definition does not coincide with that of $h$. So we will consider an other quasiconformal map $H$, quasiconformally conjugated to $\widetilde{\beta}$, thus preserving an other Beltrami form $\mu'$.

We had called $U$ the upper chessboard box, and $\phi$ was an analytic isomorphism from $\bbC\setminus \adh{U}$ to $\bbC\setminus\adh{\bbH}$, that commutes with $T_1$. The map $\widetilde{\beta}$ is equal to $T_\tau\circ T_{-v} \circ h \circ \phi^{-1}$ on $\bbC\setminus\adh{\bbH}$, and to some quasiconformal ``rotation'' on $\adh{\bbH}$. So the domain of definition of $\widetilde{\beta}$ is $\adh{\bbH}\cup\phi(\operatorname{Def}(h))$.
Remind that $h$ is a horn map whose phase we do not care about. The map $T_\tau\circ T_{-v} \circ h$ is also a horn map, and \emph{for convenience, from now on, $h$ will refer to that one}.
We want to conjugate $\widetilde{\beta}$ back by $\phi$. So we must first extend $\phi$ on $\adh{U}$ into a quasiconformal homeomorphism $\widetilde{\phi}$ of $\bbC$, and will then define
\[
 H = \widetilde{\phi}^{-1} \circ \widetilde{\beta} \circ \widetilde{\phi}
\]
and the following $H$-invariant and $T_1$-invariant Beltrami form
\[
 \mu' = \widetilde{\phi}^*(\mu)
\]
We do not have much choice on the definition $\widetilde{\phi}$, because of the requirement in propositions~\ref{prop_mucomp}+\ref{prop_pullpsih} that $\mu'$ be the pull-back by $h$ of some Beltrami form and be $T_1$-invariant.

\begin{lemma}
 There is one and only one extension $\widetilde{\phi}$ of $\phi$ such that $\mu'$ is a pull-back by $h$, and it has the form
 \[
   \widetilde{\phi}\big|_{\adh{U}} = (\widetilde{\beta}\big|_{\adh{\bbH}})^{-1} \circ h
 \]
 This is also the only extension such that $H$ is ``compatible'' with $h$, i.e. has the form $H = \operatorname{sg} \circ h$.
 Moreover, the following holds on all of $\operatorname{Def}(h)$:
 \[
  H = \widetilde{\phi}^{-1} \circ h
 \]
\end{lemma}
\varthm{Proof}
 Indeed, let $\widetilde{\phi}$ be any extension of $\phi$ such that the Beltrami form $\mu'=\widetilde{\phi}^* \mu$ is a pull-back by $h$. Note that $\mu' = \phi^* \mu $ on $\operatorname{Def}(h)\setminus U$. Now, the chessboard boxes of $h$ are mapped isomorphically by $h$ to half planes delimited by $\bbR$. There are infinitely many mapped to each one, so there is certainly one mapped to the upper half plane and different from $U$. Let us call it $V$. The form $\mu'$ is fixed on $V$, and it is pushed by the branch $h|_V$ to a fixed form on $\bbH$. The pull-back by the branch $h|_U$ is thus fixed: it means $\mu'$ is unique.
Thus any $\widetilde{\phi}$ has fixed Beltrami differential on $U$ and since it is extends $\phi$, it has fixed values on the boundary of $U$, thus it is unique. 
 
Let us now consider the compatibility requirement. Assume that a function $\widetilde{\phi}$ was defined such that the map $H$ defined by $H=\widetilde{\phi}^{-1}\circ \widetilde{\beta}\circ\widetilde{\phi}$ has the form $H=\operatorname{sg}\circ h$. Then, $\mu'=\widetilde{\phi}^*(\mu)$ being $H$-invariant, it is certainly a pull-back by $h$. Since $H = \widetilde{\phi}^{-1} \circ \widetilde{\beta} \circ \widetilde{\phi}$ on $\bbC$ and $\widetilde{\beta} = h \circ \phi^{-1}$ on $\bbC\setminus \bbH$, we have $H = \widetilde{\phi}^{-1} \circ h$ on $\bbC\setminus\bbH$. But we know that $h(\bbC\setminus\bbH)$ equals all of $\bbC$. The requirement $H=\operatorname{sg}\circ h$ on $\bbC$ thus implies that $\operatorname{sg} = \widetilde{\phi}^{-1}$ on all of $\bbC$. Thus $\widetilde{\phi}^{-1} \circ \widetilde{\beta} \circ \widetilde{\phi} = H = \widetilde{\phi}^{-1} \circ h$ on all of $\bbC$, which simplifies into $\widetilde{\beta} \circ \widetilde{\phi} = h$. Taken on $\adh{U}$, this yields the formula $\widetilde{\phi}\big|_{\adh{U}} = (\widetilde{\beta}\big|_{\adh{\bbH}})^{-1} \circ h$.

Conversely, one checks that this formula works: in other words, this $\widetilde{\phi}$ extends $\phi$ to a quasiconformal homeomorphism of $\bbC$ commuting with $T_1$, such that the relation $H = \widetilde{\phi}^{-1} \circ h$ holds on all of $\operatorname{Def}(h)$, and thus the $H$-invariant Beltrami form $\mu'$ is a pull-back by $h$.
\qed\medbreak

Let $\widetilde{\phi}$ be this extension. Let $s=S \circ \widetilde{\phi}$ be the conjugacy from $H$ to $l$. The map $s$ is quasiconformal, commutes with $T_1$, and fixes $w_0$. We have the following commutative diagram:
\[
 \xymatrix@!{
                                & \operatorname{def}(h) \ar[r]^s \ar[d]^H \ar[dl]_h & \operatorname{def}(l) \ar[d]^l \\
   \bbC \ar[r]_{\widetilde{\phi}^{-1}} & \bbC        \ar[r]_s                    & \bbC                 \\
 }
\]
Whose outer part can be read
\[
 \xymatrix@!{
  \operatorname{def}(h) \ar[r]^s \ar[d]_h            & \operatorname{def}(l) \ar[d]^l  \\
  \bbC        \ar[r]_{s\circ\widetilde{\phi}^{-1}} & \bbC                  \\
 }
\]
We now apply proposition~\ref{prop_l} to it. This proves theorem~\ref{thm_lish}.\qed

\

Thus $s(\operatorname{def}(h))= \operatorname{def}(l)= \operatorname{def}(h)$. As noted in the proof of proposition~\ref{prop_l}, $s$ satisfies the conditions of proposition~\ref{prop_mucomp} (the Beltrami form straightened by $s$ is $\mu'$). Thus $s$ fixes every point in $\psi_+^{-1}(J')$.

\begin{cor}
 The Siegel disk $\Delta$ of the map $\widehat{h}_\sigma$ is compactly contained in the domain of definition of $\widehat{h}_\sigma$. It's boundary is a quasicircle, which contains the upper main critical point of $\widehat{h}_\sigma$ (the one which is on the boundary of the uppermost chessboard box), and no other. 
\end{cor}

Let us also use the notation $\Delta$ for the preimage of $\Delta$ by the covering map $\bbC \to \bbC/\bbZ$. All chessboard boxes $V$ are mapped by $h$ to the upper or the lower half plane. Let us call $V$ \emph{light} in the first case, and \emph{dark} in the second case.

\begin{cor}[combinatorial description]
 The image of the light chessboard boxes of $\operatorname{def}(h)$ by the quasiconformal isomorphism $s$ gives the connected components of the first preimage by $h_\sigma$ of the Siegel disk $\Delta$. From each of these components, $h_\sigma$ is an isomorphism to $\Delta$. All the critical points of $h$ belong to the image by $s$ of the chessboard graph (the analytic lines separating the boxes). The dark boxes map by $s$ to components mapped isomorphically to $\bbC\setminus\adh{\Delta}$ by $h_\sigma$.
\end{cor}

So far, we analyzed the dynamics of the horn map $h_\sigma$, which live in the Fatou coordinates (more precisely the universal cover of the Ecalle-Voronin cylinders). Let us give a informal description of what happens in the dynamical plane where $P$ lives. A more precise version of this analysis is given in my thesis \cite{C}. We saw that $\mu'$ is the pull-back by $\psi_+$ of some Beltrami form. Let us call it $\mu''$. 
Let us have a closer look at the complement~\ref{comp_R} of proposition~\ref{prop_mucomp}.
The form $\mu''$ has a straightening $R$ preserving $P$, every point of $J$, every component of $\bbC\setminus J$, and every point in the grand orbit of the critical point. The map $R$ satisfies the relation $R\circ\psi_+ = \psi_+ \circ s$. It sends $\Delta$ to what is called the \emph{Virtual Siegel disk} $\Delta'$ of the (modified) Lavaurs map $g_\sigma$. It is periodic of period $q$ under $P$, and fixed by $g_\sigma$. We recall that $g_\sigma : \insb{K} \to \bbC$ is holomorphic and satisfies the relation $g_\sigma \circ \psi_+ = \psi_+ \circ h_\sigma$. The filled-in Julia-Lavaurs set $K(P,g_\sigma)$ is by definition the complement of the escaping set, the latter being the union for $n\in\bbN$ of $g_\sigma^{-n}(\bbC\setminus K)$. The Julia-Lavaurs set $J(P,g_\sigma)$ is by definition the boundary of $K(P,g_\sigma)$ in $\bbC$. See \cite{D} and \cite{L}. The map $R$ maps the chessboard graph of $\insb{K}$ into a subset of the Julia-Lavaurs set $J(P,g_\sigma)$, and every light box to a component of the first preimage of $\Delta'$ by $g_\sigma$, or equivalently to an iterated preimage by $P$ of $\Delta'$. The following pictures illustrate all this.

\subsection{Illustrations}

We illustrate below the case of the Lavaurs maps.

\begin{figure}
 \begin{center}
 \begin{picture}(160,460)
  \rotatebox{90}{%
  \put(45,0){%
   \put(-45,145){%
    \scalebox{0.68}{\rotatebox{270}{%
     \includegraphics*[20,0][230,307]{Pics/r_Rab_5_Ech_Fatou.eps}
   }}}
   \put(195,2){%
    \includegraphics*[0,10][210,154]{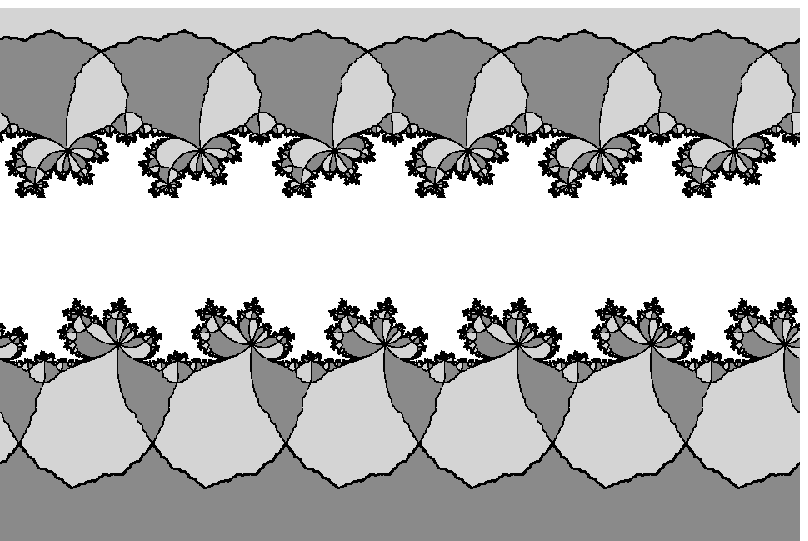} 
   }
   \put(169,70){\scalebox{1.3}{$\overset{s}{\longrightarrow}$}}
   \put(298,150){$\Delta$}
  }
  }
 \end{picture}
 \end{center}
 \caption{The parabolic chessboard in Fatou coordinates, and it's
   image by the straightening $s$, for $p/q = 2/5$, and
   $\displaystyle\theta= \frac{\sqrt{5}-1}{2}$ (the golden mean). The
   Siegel disk is at the upper end. The figure has been rotated by
   $90$~degrees.}
\end{figure}

\begin{figure}
 \begin{center}
 \begin{picture}(274,450)
  \rotatebox{90}{%
  \put(0,10){%
   \put(0,0){%
    \includegraphics{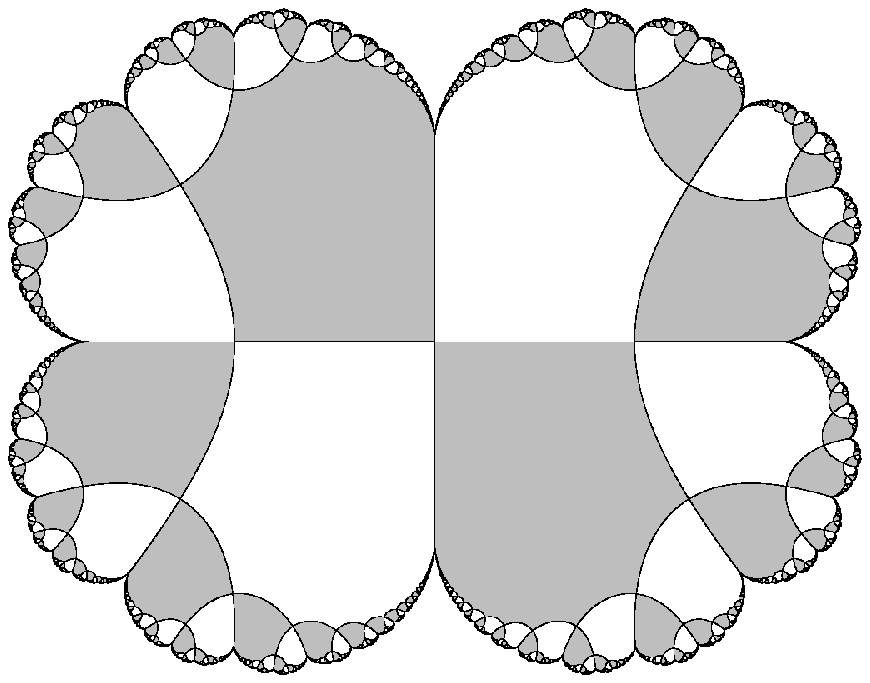} 
   }
   \put(220,0){%
    \includegraphics*[0,27][230,280]{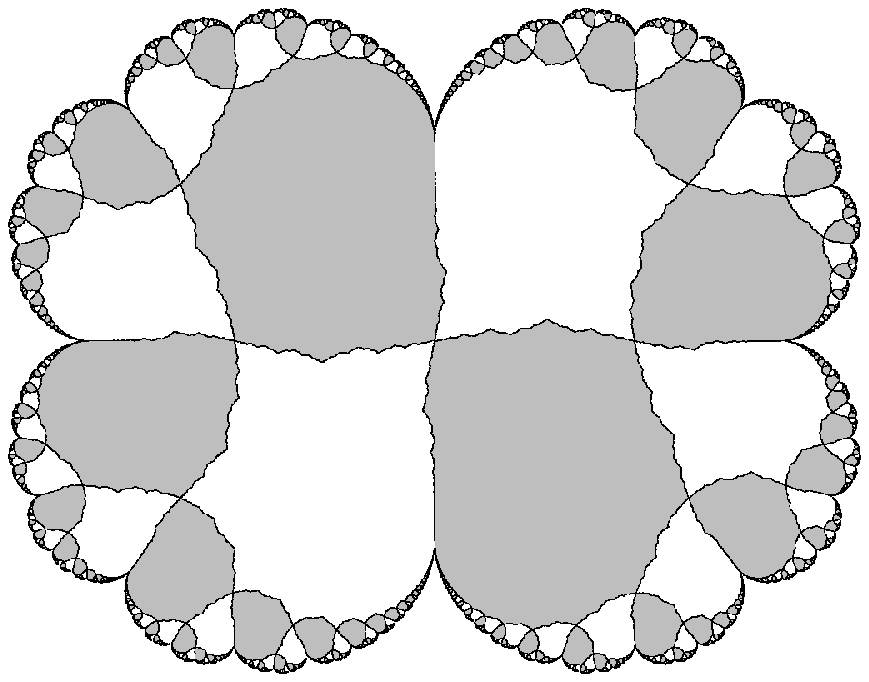} 
   }
   \put(210,122){\scalebox{1.3}{$\overset{R}{\longrightarrow}$}}
   \put(370,152){$\Delta'$}
  }
  }
 \end{picture}
 \end{center}
 \caption{The parabolic chessboard in initial coordinates, for $p/q =
   0/1$, and it's image by the straightening $R$, with $\theta=$ the
   golden mean. The figure has been rotated by $90$~degrees.}
\end{figure}

\clearpage

\

\begin{figure}[h]
 \begin{center}
  \includegraphics*[0,17][230,290]{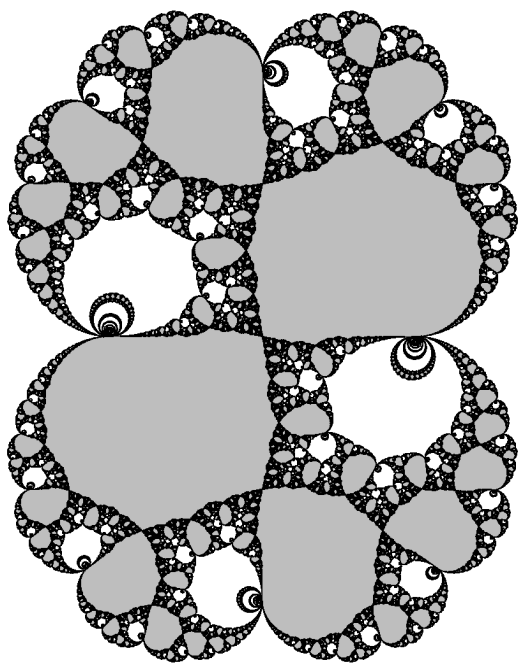} 
 \end{center}
 \caption{The corresponding filled-in Julia-Lavaurs set. It is the non-escaping set under $P$ and $g_\sigma$ (with the appropriate value of $\sigma$. The picture shows the boundary in black, and the interior in gray.}
\end{figure}

\pagebreak

\begin{figure}
 \begin{center}
 \includegraphics{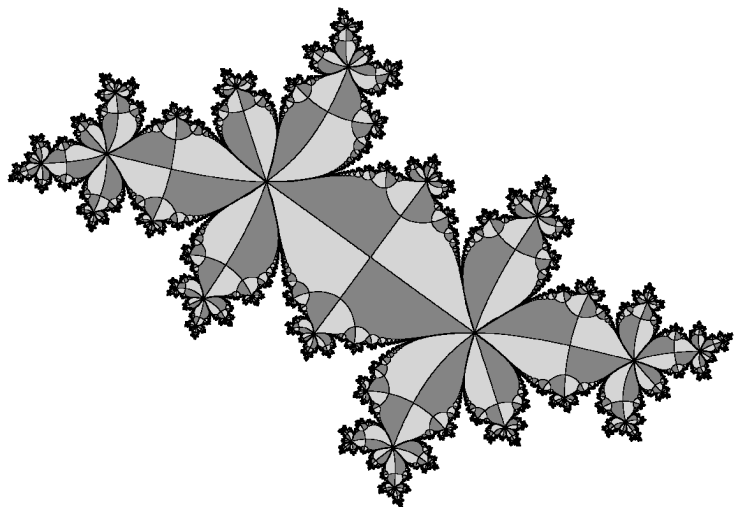}
 \end{center}
 \caption{The Parabolic chessboard of the Julia set in initial coordinates, for $p/q=2/5$ and $\theta=$ the golden mean.}
\end{figure}

\begin{figure}
 \begin{center}
  \includegraphics{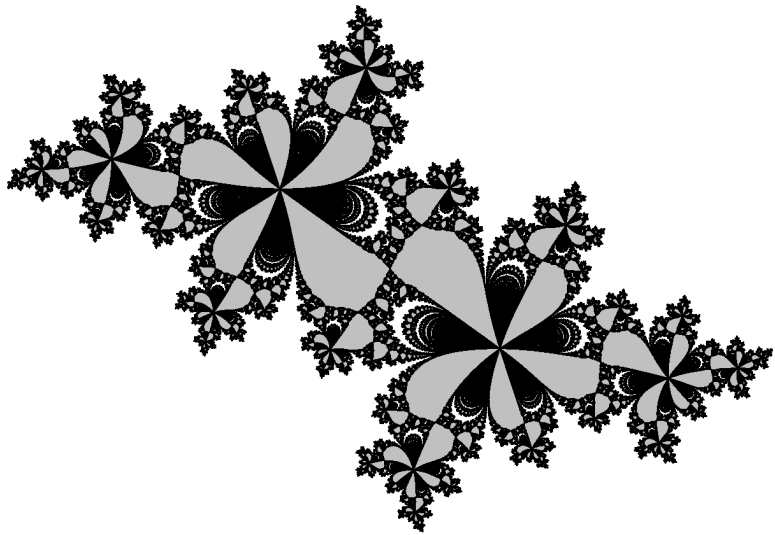}
 \end{center}
 \caption{There follows the Julia-Lavaurs set for $p/q=2/5$ and
   $\theta=$ the golden mean.}
\end{figure}

\clearpage


\section{Getting more control for Lavaurs maps}

The notion of $c$-quasisymmetrical function which is used here is not invariant by M\"obius transformations of $\bbS^1$ for a given value of $c>1$: by definition an orientation preserving homeomorphism $f: \bbT \to \bbT$ is $c$-quasisymmetric if and only if, calling $\mathbf{f}$ a lift of $f$ to a homeomorphism of $\bbR$, we have
\[
 \forall x\in\bbR,\ \forall h\in]0,1[,\ c^{-1} \leq \frac{\mathbf{f}(x+h)-\mathbf{f}(x)}{\mathbf{f}(x)-\mathbf{f}(x-h)} \leq c
\]

\begin{prop}[uniformity in the Herman, \'Swi\c{a}tek theorem]\label{prop_unifHeSw}
 In theorem~\ref{thm_HeSw}, the constant $c$ only depends on the following data:
\[M,W,l,A,A'\text{ and }M'\]
\end{prop}
I claim that this statement can be deduced from \cite{Pe2}.

\begin{prop}[uniformity in the Ahlfors, Beurling theorem]
 There exists a function $K(k) : ]1,+\infty[ \to ]1,+\infty[$, such that $K(k) \tend 1$ when $k \tend 1$, and such that every $k$-quasisymmetrical homeomorphism $f$ of the unit circle extends to a homeomorphism $F$ from $\adh{\bbD}$ to $\adh{\bbD}$, $K$-quasiconformal in $\bbD$, and that moreover fixes $0$.
\end{prop}
\varthm{Proof}
 The statement is in \cite{AB}, except the "fixing $0$" claim. For any given $k$ and $K$, the class of $K$-quasiconformal homeomorphisms of $\bbD$ that extends to a $k$-quasisymmetrical homeomorphism of the unit circle, is compact for the uniform convergence on $\bbD$. In particular, the image of $0$ by such a homeomorphism belongs to some compact $C$ of $\bbD$, which depends only on $k$ and $K$. It is then possible to bring it back to $0$ by post-composing with a $C^\infty$-diffeomorphism of $\bbD$ that is identity on a neighborhood $V$ of $\partial\bbD$, and whose ellipticity is bounded by $K'$, where $V$ and $K'$ only depend on $C$.
\qed\medbreak

We need now to introduce the following Blaschke product:
\[
 \mcP(z)=\frac{3z^2+1}{z^2+3}
\]
It has degree $2$. It leaves $\bbD$, $\bbS^1$ and $\bbS^2\setminus \adh{\bbD}$ invariant. It restricts to a degree $2$ self-covering of $\bbS^1$. It has a unique fixed point $z=1$, which is parabolic with multiplier $1$ and which has $2$ attracting axes, $[1,+\infty[$ and $]-\infty,1]$. The basin of attraction of the left axis is $\bbD$ and the basin of the right axis is $\bbS^2\setminus \adh{\bbD}$. Thus, the Julia set is $\bbS^1$.

It has a horn map that we will denote $\hbl$. It is defined on $\bbC\setminus\bbR$ and maps to all of $\bbC$. The following picture illustrates the parabolic chessboard of this map on the half plane $\bbH$. It is followed by the chessboard of a connected component of domain of definition of the horn map $h$ of the polynomial $P_{p/q}$ (for $p/q=2/5$). They are linked by the following proposition.

\begin{center}
 \begin{picture}(192,230)
  \put(-10,0){%
   \put(0,150){%
    \includegraphics{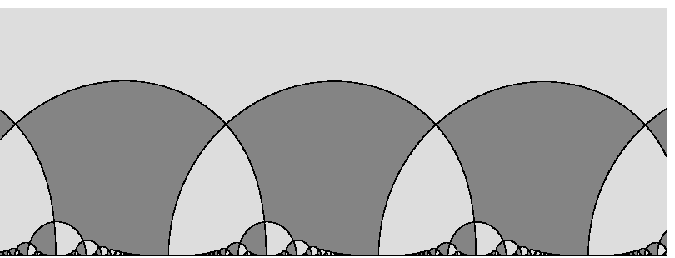}
   }
   \put(0,0){%
    \includegraphics{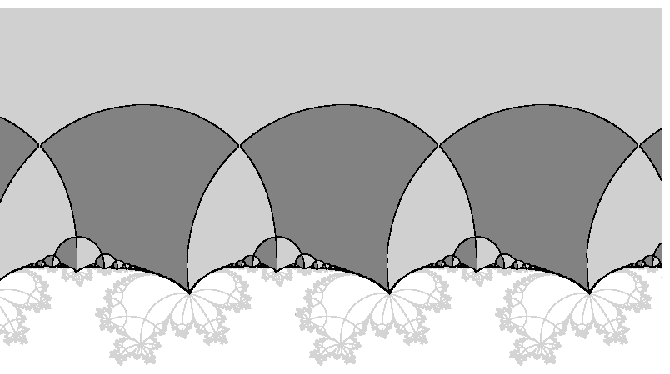}
   }
   \put(88,145){%
    \xymatrix{ \ \\ \ \ar[u]_{\xi} }
   }
  }
 \end{picture}
\end{center}

\begin{prop}[Universality of horn maps as coverings]\label{prop_univ}
 For all $p/q$, there exists a $\sigma_0$ depending on the many choices in Fatou coordinates, and an isomorphism $\xi$ from the component $U$ of $\psi_+^{-1} \ins{K}(P)$ that contains an upper half plane to the upper half plane $\bbH$, that commutes with $T_1$ and such that 
 \[
  h_{\sigma_0} = \hbl \circ \xi
 \]
 holds on $U$.
\end{prop}
A proof of this claim can be found in~\cite{C}. Basically, it is a transposition to the horn maps of the following universality (see~\cite{O}): for any quadratic polynomial $Q$ having a parabolic point, for any periodic component $C$ of $\ins{K}(Q)$, let $n$ be the minimal period. Then there is an analytic isomorphism from $C$ to $\bbD$ conjugating $Q^n$ to $\mcP$.

\

In particular, proposition~\ref{prop_univ} implies that, modulo $\bbZ$, the connected component of the domain of definition of the pre-model map $\beta$ that contains the real axis (this is a $\bbR$-symmetric annulus $A_{p/q}$), has modulus equal to $2m$ where $m$ is the modulus between the upper chessboard box and the boundary of the domain of definition of $h$, which is the same as the analog modulus for $\hbl$, thus the modulus of $A_{p/q}$ does not depend on $p/q$: it is universal.
Let us call $\bbl$ the pre-model map constructed from the Blaschke product's horn map $\hbl$, let $A$ it's central annulus of definition ($\bmod\ \bbZ$), and let $\phibl$ be to $\hbl$ what $\phi$ is to $h$.
Then, the map $\phibl \circ \xi \circ \phi^{-1}$ has a Schwarz reflection extension across $\bbR$, that we will call $\zeta$, and which is an isomorphism from the annulus $A_{p/q}$ to $A$.
We have 
\[
  \beta\big|_{A_{p/q}} = \bbl \circ \zeta
\]

From this, one deduces that, for a fixed bound $M$ on the entries of the continued fraction expansion of $\theta$, the conditions of proposition~\ref{prop_unifHeSw} are uniformly satisfied, i.e.\ with constants that do not depend on $p/q$.
We have thus proved the following theorem.

\begin{theo}[Universal bound on dilatation ratio]\label{prop_modindep}
 For all $M>0$, there exists $K>0$ such that for all $p/q$ irreducible and all irrational rotation number $\theta$ of type bounded by $M$, the straightenings $S$, $R$ and $s$ have dilatation ratio bounded by $K$.
\end{theo}

Here is an example of consequence:

\begin{cor}
 The modulus of the annulus in the cylinder $\bbC/\bbZ$ separating $\Delta$ from the boundary of the domain of definition of $h$, is bounded below (and above), independently of $p/q$, by a constant depending only on $M$.
\end{cor}
\varthm{Proof}
 This annulus is the image by $S$ of the lower half of the annulus $A_{p/q}$, and we saw this half has universal modulus.
\qed\medbreak


\section{Examples}

\subsection{Uncountably many inequivalent entire maps}\label{subsec_ex_entire}

We will provide here an uncountable set of entire maps that are topologically inequivalent, and satisfy the conditions of section~\ref{subsec_const_entire}. They will in fact satisfy the stronger following properties: $0$ and $1$ are critical values, they are the only ones, there are no asymptotic values, $0$ has at least one non critical preimage $z$ such that the component of the preimage of the unit disk that contains $z$ is bounded.

It is enough to construct inequivalent quasiregular maps $f : \bbC \to \bbC$. Indeed, by letting $\mu$ being the pull-back by $f$ of the null Beltrami form, and $\phi$ a straightening ($\phi_*(\mu) = 0$), we obtain an entire map $g = f \circ \phi^{-1}$, equivalent to $f$. 

Let $a_n \in \{0,1\}$ be indexed by $n \in \bbN^*$. For reasons that will be explained below, we require that $a_n =1$ infinitely many times. For $n\in \bbZ$, let $S_n$ be the strip defined by ``$\Re(z)\in [n,n+1]$''. We will define below two continuous maps $A$ and $B$ from $S_0$ to $\bbC$, quasiregular in the interior, and such that $\forall y\in\bbR$, $A(i\, y)=A(1+i\, y)=B(i\, y)=B(1+i\, y)$. Then $f(z)$ will be defined by $\forall n\in\bbZ$ and $\forall z\in S_n$, $f(z) = A(z-n)$ or $f(z)=B(z-n)$, according to the following rule:
\begin{itemize}
 \item if $n<0$, then we take $A$
 \item if $n=0$, then we take $B$
 \item if $n>0$, then we take $A$ if $a_n=0$, and $B$ if $a_n=1$
\end{itemize}

\clearpage 

\

\begin{figure}[h]
 \begin{center}
 \rotatebox{90}{
 \begin{picture}(450,180)
  \put(0,0){%
   \put(-2,0){%
    \put(30,165){$D_{-2}$}
    \put(105,165){$D_{-1}$}
    \put(180,165){$D_{0}$}
    \put(255,165){$D_{1}$}
    \put(330,165){$D_{2}$}
    \put(405,165){$D_{3}$}
   }
   \put(-2,0){%
    \put(30,5){$A$}
    \put(105,5){$A$}
    \put(180,5){$B$}
    \put(255,5){$B$}
    \put(330,5){$A$}
    \put(405,5){$B$}
   }
   \put(0,20){%
    \put(0,0){%
     \includegraphics{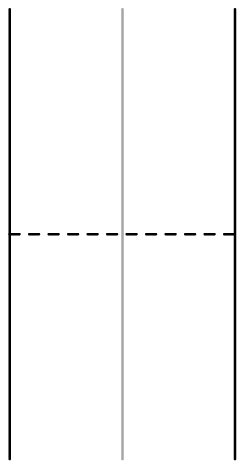}
    }
    \put(75,0){%
     \includegraphics{Pics/A.ps}
    }
    \put(150,0){%
     \includegraphics{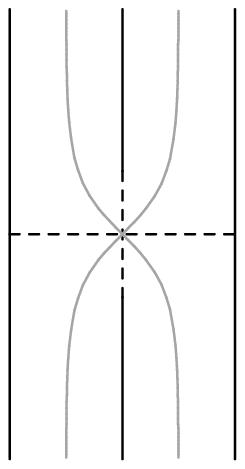}
    }
    \put(225,0){%
     \includegraphics{Pics/B.ps}
    }
    \put(300,0){%
     \includegraphics{Pics/A.ps}
    }
    \put(375,0){%
     \includegraphics{Pics/B.ps}
    }
   }
  }
 \end{picture}
 }
 \end{center}
 \caption{Example of sequence for the map $f$. The figure has been
   rotated by $90$~degrees.}
\end{figure}

\pagebreak

Let us now define these maps $A$ and $B$.
\[B(z) = 1-(\cos(\pi z))^4\]
Note that $B$ is holomorphic. It is illustrated by the following picture.
\begin{center}
 \scalebox{0.8}{%
 \begin{picture}(350,140)
  \put(-45,0){
  \scalebox{1}{
   \put(0,0){
    \includegraphics{Pics/B.ps}
   }
  }
  \put(72,64){$\overset{-\cos (\pi z)}{\longrightarrow}$}
  \put(110,0){%
   \scalebox{0.7}{
    \fbox{
     \includegraphics{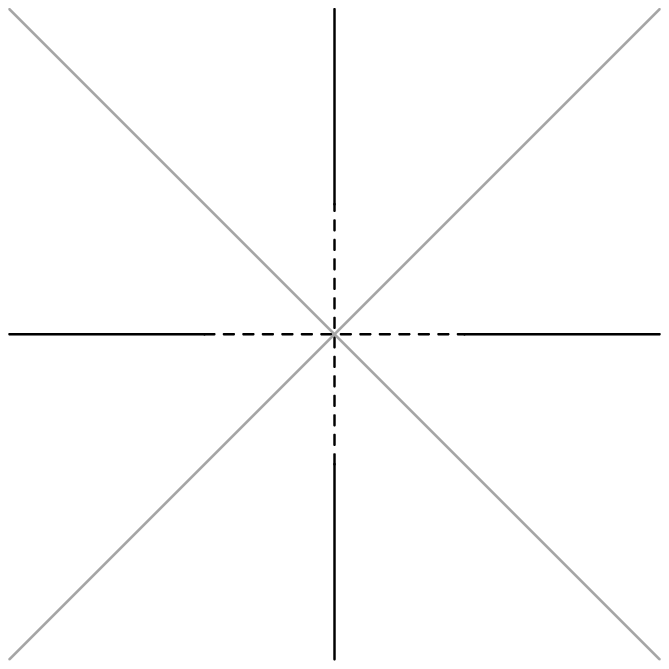}
    }
   }
  }
  \put(262,64){$\overset{1-z^4}{\longrightarrow}$}
  \put(290,0){%
   \scalebox{0.86}{
    \fbox{
     \includegraphics[231,320][380,470]{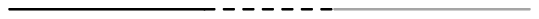}
    }
   }
  }
  \put(347,55){$0$}
  \put(377,55){$1$}
  }
 \end{picture}
 }
\end{center}
Let
\[A(z) = \frac{1-\cos(2\pi h(z))}{2}\]
where
\[h(x+i\, y) = x+i\, l(y)\]
and $l$ is a increasing homeomorphism of $\bbR$ such that $A$ and $B$ coincide on $i \bbR$. An explicit computation yields
\[1-(\cosh(\pi y))^4 = \frac{1-\cosh(2\pi l(y))}{2}\]
from which one deduces that $l$ is a diffeomorphism, with derivative satisfying $l'(z)\in [\sqrt{2},2[$.
Thus $A$ is $K$-quasiconformal with $K = 2$.
\begin{center}
 \begin{picture}(250,140)
  \scalebox{1}{
   \put(0,0){%
    \includegraphics{Pics/A.ps}
   }
  }
  \put(85,64){$\overset{A}{\longrightarrow}$}
  \put(120,0){%
   \put(0,0){%
   \scalebox{0.86}{
    \fbox{
     \includegraphics[231,320][380,470]{Pics/p0.eps}
    }
   }
   }
   \put(57,55){$0$}
   \put(87,55){$1$}
  }
 \end{picture}
\end{center}

Now let us prove that two different sequences $a_n$ yield two inequivalent maps. Let us consider any self-avoiding path $\gamma$ from $0$ to $1$, in $\bbC\setminus \{0,1\}$. Any two such paths are isotopic. As a consequence, starting from any preimage $z$ of $0$, resp.\ 1, let us say that the local degree is $d$, the set of points one attains by following one of the $d$ lifts of $\gamma$, resp.\ $\gamma^{-1}$, starting from $z$ do not depend on $\gamma$, and is thus covariant by topological equivalences.
This gives a graph embedded in the complex plane, with vertices labeled $0$ or $1$, that we will call the skeleton. 
\begin{center}
 \xymatrix@!@=10pt{
   & & & & & 0 & & 0 & & & & 0 & \\
   \ldots \ar@{-}[r] &
   1 \ar@{-}[r] & 
   0 \ar@{-}[r] & 
   1 \ar@{-}[r] &
   0 \ar@{-}[r] &
   1 \ar@{-}[r] \ar@{-}[u] \ar@{-}[d] &
   0 \ar@{-}[r] &
   1 \ar@{-}[r] \ar@{-}[u] \ar@{-}[d] &
   0 \ar@{-}[r] &
   1 \ar@{-}[r] &
   0 \ar@{-}[r] &
   1 \ar@{-}[r] \ar@{-}[u] \ar@{-}[d] &
   \ldots \\
   & & & & & 0 & & 0 & & & & 0 & \\
 }
\end{center}
This graph has a spine, a bi-infinite sequence of vertices alternately labeled $0$ and $1$, each linked to the other. All the $1$s are on the spine, let us call them the vertebrae. Some vertebrae have ribs, i.e.\ edges that end with lone vertices labeled $0$: these are the vertebrae associated to the map $B$. The spine has no ribs before the one associated to $n=0$, it has one there, and infinitely many ribs after, since we required that $a_n=1$ for infinitely many $n>0$.  Equivalent maps must have homeomorphic skeleton, with homeomorphism respecting labeled points (in fact, even better than that). But since it is possible to topologically identify which $n\in \bbZ$ is associated to each vertebra, this implies equivalent maps have the same sequence $a_n$.

The conclusion follows, because there are uncountably many sequences $a \in \{0,1\}^{\bbN^*}$ with infinitely many $1$s.

Remark: in fact, I had first thought of the entire maps $g$ by gluing pieces of ramified coverings, but to prove that the Riemann surface thus constructed is isomorphic to $\bbC$ and not to $\bbD$, it is enough to prove it is quasiconformally equivalent to $\bbC$, which amounts to directly define a quasiregular map $f$.

\

More generally, simply connected ramified covers over $\bbC$, ramified only over $0$ and $1$, are completely classified by their skeleton (up to a positive homeomorphism of the plane in which the skeleton is drawn), which can be any finite or infinite tree in the plane (connected graph with no cycles), with no accumulation, and with at least one edge. Rigorously, you need to label the vertices alternately by $0$ and $1$, but note that any unlabeled skeleton has exactly two labelings. Now, when is the Riemann surface thus constructed isomorphic to $\bbC$ instead of $\bbD$~? Is it possible to find a simple criterion to be read on the skeleton~?
Note that these skeletons are very convenient for understanding polynomials. Their number of edges is equal to the degree.


\clearpage

\section{Objects of the parabolic implosion}\label{sec_obj}

Here, we give the definition of the objects involved in parabolic implosion.

\

\begin{center}
 \scalebox{0.6}{%
  \includegraphics{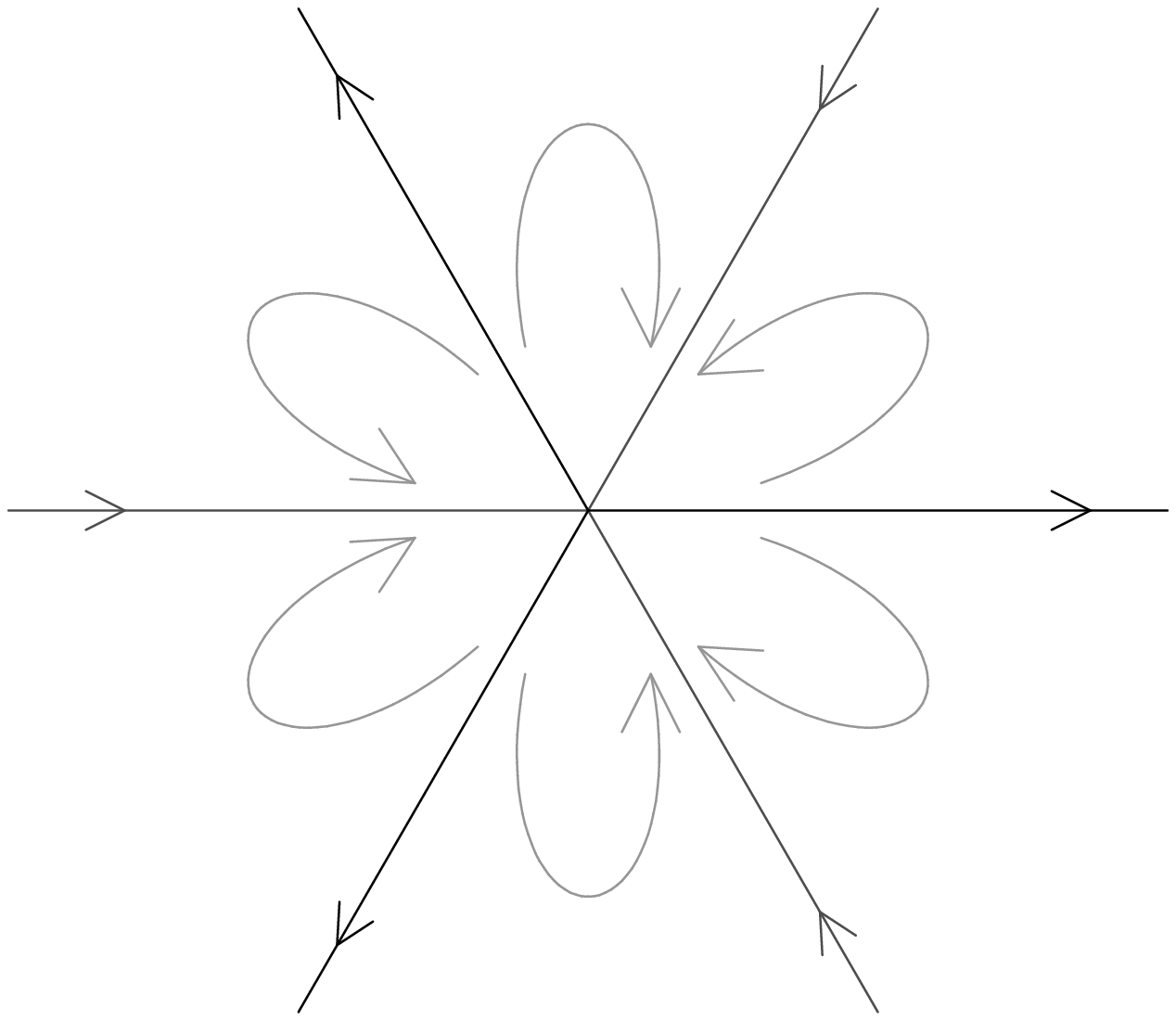}%
 }
\end{center}

\

The repelling axes of the parabolic point $z=0$ are labeled by $j$, which can take $q$ values.
The attracting ones are labeled by $i$, which can also take $q$ values.
These $2q$ axes are half lines starting from $0$, alternatively attracting and repelling, and they all make equal angles.
The \emph{extended repelling Fatou parameterization} associated to the axis $j$ is noted
\[
 \psi_{+,j} : \bbC_j \to \bbC
\]
where $\bbC_j$ is just a copy of $\bbC$ that we labeled $j$. It is a holomorphic map, semi-conjugating $T_1$ to $P^q$: the following diagram commutes
\[
 \xymatrix{
  \bbC_j \ar[d]_{\psi_{+,j}} \ar[r]^{T_1} & \bbC_j \ar[d]^{\psi_{+,j}}\\
  \bbC \ar[r]_{P^q} & \bbC \\
 }
\]
The extended attracting Fatou coordinate will be noted
\[
  \phi_{-,i} : \ins{K}_i \to \bbC_i
\]
where $\ins{K}_i$ stands for the set of points in $\bbC$ attracted by the axis $i$ under the dynamics of $P^q$ (note the iterate). The union of the $\ins{K}_i$ is disjoint, and is equal to $\ins{K}$. We have the following commuting diagram
\[
 \xymatrix{
  \bbC_i  \ar[r]^{T_1} & \bbC_i \\
  \ins{K}_i \ar[u]^{\phi_{-,i}} \ar[r]_{P^q} & \ins{K}_i \ar[u]_{\phi_{-,i}}\\
 }
\]

Moreover, $P$ acts on the set of repelling (resp.\ attracting) axes, like a rotation of angle $p/q$ turns. Let us note $i^P$ the axis $i$ maps to, and $j^P$ the axis $j$ maps to. We have the following properties:
\[
 \xymatrix{
  \bbC_i \ar[r]^{T_{1/q}} & \bbC_{i^P} & &
  \bbC_j \ar[d]_{\psi_{+,j}} \ar[r]^{T_{1/q}} & \bbC_{j^P} \ar[d]^{\psi_{+,j^P}} \\
  \ins{K}_i \ar[u]^{\phi_{-,i}} \ar[r]_{P} & \ins{K}_{i^P} \ar[u]_{\phi_{-,i^P}} & &
  \bbC \ar[r]_{P} & \bbC \\
 } 
\]

With these properties, the Fatou coordinates and parameterizations are all unique up to a composition with two translations: more precisely, the $\psi_{+,j}$ may be all pre-composed with the same translation, and the $\phi_{-,i}$ may be all post-composed with an other common translation. (Thus there are two $\bbC$-degrees of freedom in the choices.)

We introduce the maps
\[
  \phi_{\div,i} : \ins{K} \to \bbC_{i}
\]
by the formula
\[
 \forall i,\ \forall i',\ \forall z\in\ins{K}_{i'},\ \phi_{\div,i}(z) = \phi_{-,i}\circ P^k (z) = T_{\frac{k}{q}} \circ \phi_{-,i'}
\]
where $k\geq 0$ is the smallest number of time $P$ must be iterated to map the axis $i'$ to the axis $i$: $k$ depends on $i$ and $i'$. (Note that the definition is slightly different than the one in \cite{C}, since here our map $\phi_{\div,i}$ is indexed by $i$ instead of $j$, and we do not include the translation $T_\sigma$.)

Let us mention the following relation:
\[
 \forall i, \forall z\in\ins{K},\ \phi_{\div,i}\circ P(z) = \phi_{\div,i}(z) + m_{i,z}
\]
where, if we call $i'$ the index such that $z\in\ins{K}_{i'}$, we define $m_{i,z}=1$ if $i'=i$ and $m_{i,z}=0$ else.

In parabolic implosion, one must choose $\nu \in \{-1,1\}$, which stands for the \emph{direction} of the implosion, positive or negative. If $\nu=1$, we define $j(i)$ to be the repelling axis right after the attracting axis labeled $i$ in the cyclic order. If $\nu=-1$, we take for $j(i)$ the one just before. The inverse map will be noted $i(j)$. The choice in the direction of the implosion is not to be mistaken with the choice of the end of the cylinder where we put the Siegel disk. 

To simplify notations, we choose an index $j=j_0$, fix it, and we call $\bbC_+ = \bbC_{j_0}$ and $\bbC_- = \bbC_{i(j_0)}$. The notation $\psi_{+}$ will stand for $\psi_{+,j_0}$, and $\phi_{\div}$ for $\phi_{\div,i(j_0)}$.
Let us consider the following \emph{non commuting} diagram
\[
 \xymatrix@=30pt{
  \save[]+<25pt,0pt>*+{\bbC}="targ" \ar[dr]^{\phi_{\div}} \restore
  & \\
  \bbC_+ \ar"targ"^{\psi_+} & \ar[l]^{T_\sigma} \bbC_{-} \\
 }
\]
where $\sigma$ is a parameter: $\sigma \in \bbC$.

The \emph{modified} Lavaurs maps are noted
\[
 g_\sigma : \ins{K} \to \bbC
\]
They consists in turning around the triangular diagram from the vertex labeled $\bbC$ to itself. In other words, they are defined by the formula
\[
 g_\sigma = \psi_+ \circ T_\sigma \circ \phi_\div
\]

The extended horn maps consist in turning around the diagram, this time starting from $\bbC_+$. In other words, 
\[
  h_{\sigma} = T_\sigma \circ \phi_\div \circ \psi_+
\]
Thus they are defined on the subset $\phi_{+}^{-1}(\ins{K})$ of $\bbC_+$:
\[
  h_{\sigma} : \phi_+^{-1}(\ins{K}) \to \bbC_+
\]
The maps $g$ and $h$, for the same value of $\sigma$, are semi-conjugated one of the other, as can be derived from either the triangle or the formulas: $g_\sigma \circ \psi_+ = \psi_+ \circ h_\sigma$ and $h_\sigma \circ (T_\sigma \circ \phi_\div) = (T_\sigma \circ \phi_\div) \circ g_\sigma$.

\varthm{Remark}
The definition of the non-modified Lavaurs map is just slightly more complicated.
\[
 \forall i,\ \forall z\in\ins{K}_i,\ L_\sigma (z) =
 \psi_{+,j(i)} \circ T_\sigma \circ \phi_{-,i}(z)
\]
This definition does not single out any particular value of $j$.
Now we have
\[
 \forall i,\ \forall z\in\ins{K}_i,\ g_{\sigma} (z) = P^k \circ L_\sigma (z)
\]
where $k$ is the smallest non negative iterate of $P$ that maps the axis $i$ to $i(j_0)$.
\medbreak


\section{Conclusion}

The model should enable proving local connexity and zero measure of the Julia-Lavaurs set of quadratic polynomials, and also help study the case of some entire maps.

The construction we propose can be carried-out for much more general maps, to give other pre-models. But we are not always guaranteed that the straightening gives the same map we started from. It would be interesting to study the case where there is less rigidity, and/or more singular values.


\section*{Aknowledgments}

This work, except for the entire maps, is part of my PhD thesis. I
would like to thank my advisor Adrien Douady. I would like to thank
Nuria Fagella, Marguerite Flexor, Lukas Geyer and Dierk Schleicher for helpful
discussions. I would also like to thank Marguerite Flexor for the
uncountably many times she read different preliminary versions of this
work.


\end{document}